\documentclass[11pt,twoside,a4paper]{article}

\usepackage{amsmath,bbm,amssymb,mathtools,geometry, xcolor,mathabx,mathrsfs,amsthm,resmes}
\usepackage{hyperref}

\allowdisplaybreaks
\mathtoolsset{showonlyrefs}

\newcommand{\D}{\nabla}
\newcommand{\dd}{\partial}
\newcommand{\norm}[1]{\left\lVert#1\right\rVert}

\newtheorem{Thm}{Theorem}[section]
\newtheorem{Lem}[Thm]{Lemma}
\newtheorem{Rem}[Thm]{Remark}

\newtheorem{Def}[Thm]{Definition}
\newtheorem{Prop}[Thm]{Proposition}

\DeclarePairedDelimiter\abs{\lvert}{\rvert}

\def\Om{\Omega}

\def\R{\mathbb R}

\def\N{\mathbb N}
\def\e{\varepsilon}
\def\ecc{\mathbf{e}^C}

\def\Xint#1{\mathchoice
   {\XXint\displaystyle\textstyle{#1}}
   {\XXint\textstyle\scriptstyle{#1}}
   {\XXint\scriptstyle\scriptscriptstyle{#1}}
   {\XXint\scriptscriptstyle\scriptscriptstyle{#1}}
   \!\int}
\def\XXint#1#2#3{{\setbox0=\hbox{$#1{#2#3}{\int}$}
     \vcenter{\hbox{$#2#3$}}\kern-.5\wd0}}

\def\dashint{\Xint-}
\newcommand{\prodscal}[2]{\left\langle{#1},{#2}\right\rangle}

\title{Regularity results for almost-minimizers of anisotropic free interface problem with H\"{o}lder dependence on the position}

\author{L. Esposito, L. Lamberti, G. Pisante}

\begin{document}
\maketitle

\begin{abstract}
 \noindent 
We establish regularity results for almost-minimizers of a class of variational problems involving both bulk and interface energies. The bulk energy is of Dirichlet type. The surface energy exhibits anisotropic behaviour and is defined by means of an ellipsoidal density that is H\"older continuous with respect to the position variable.
\end{abstract}


	
\makeatletter

\makeatother

\section{Introduction and statements}
The existence and regularity of solutions to variational problems, encompassing both bulk and interface energies, have been extensively studied across various disciplines and remain a focal point of much mathematical research. These problems serve to describe a broad spectrum of phenomena in applied sciences, including nonlinear elasticity, materials science and image segmentation in computer vision (see for instance \cite{AC,AR,CNT,FF,FFLM,FJ,Gur,JP,Lar}).
A model integral functional initially introduced to study minimal energy configurations of two conducting materials by R. V. Kohn $\&$ G. Strang and F. Murat $\&$ L. Tartar, \cite{KS, MT}, and later recovered by L. Ambrosio $\&$ G. Buttazzo, and F. H. Lin, \cite{AB, Lin}, is the following:
\begin{equation}\label{Mod1}
\int_\Om\sigma_E(x)|\nabla u|^2\,dx+P(E;\Om),
\end{equation}
where $\sigma_E:=\alpha \mathbbm{1}_E+\beta \mathbbm{1}_{\Omega\setminus E}$, $0<\alpha<\beta $,  with $E\subset\Omega\subset \R^n$ and $u\in H^1(\Omega)$. Here, $\mathbbm{1}_E$ stands for the characteristic function of $E$ and $P(E;\Om)$ denotes the perimeter of the set $E$ in $\Om$.
In \cite{AB, Lin},  the authors proved the existence and  regularity for   minimal configurations $(E,u)$ of \eqref{Mod1}. In a broader context, F. H. Lin $\&$ R. V. Kohn addressed more generalized Dirichlet energies as outlined in \cite{LK},
\begin{equation}\label{Mod2}
\int_\Om \left(F(x,u,\nabla u)+\mathbbm{1}_{E}G(x,u,\nabla u)\right)\,dx+\int_{\Omega\cap\dd^{\ast} E} \Psi(x,\nu_E(x))\;d\mathcal{H}^{n-1}(x),
\end{equation}
 with the  constraints
$$u=u_0\,\,\,\mathrm{on}\,\,\,\partial\Omega\,\,\,\mathrm{and}\,\,\, |E|=d.$$
Here, $\nu_E$ is the measure-theoretic outward unit normal to the reduced boundary $\dd^*E$ of $E$.
The regularity for the minimizing pair $(E,u)$ is quite intricate to establish, especially concerning the free boundary $\partial E$ due to the interaction between the bulk term and the perimeter term. In order to illustrate the regularity results of the free boundary, we introduce the following notations.
We define the set of regular points of $\partial E$ as follows:
\begin{equation}
\mbox{Reg}(E):= \left\{x\in \partial E\,:\,\partial E \text{ is a }C^{1,\gamma} \text{ hypersurface in } B_{\varepsilon}(x)\subset \Omega,\text{ for some } \varepsilon>0 \text{ and }\gamma\in(0,1)
\right\},
\end{equation}
and accordingly we define the set of singular points of $\partial E$ as
\begin{equation}
\Sigma(E) := (\partial E \cap \Omega)\setminus \mbox{Reg}(E).
\end{equation}
The most notable advancements in regularity results regarding the free boundary $\partial E$ of minimizers of the functional \eqref{Mod1} have been accomplished by G. De Philippis $\&$ A. Figalli in \cite{DF} and N. Fusco $\&$ V. Julin in \cite{FJ}.
They proved that for minimal configurations of the functional \eqref{Mod1} it turns out that
\begin{equation}
\label{redu}
\operatorname{dim}_{\mathcal{H}}(\Sigma(E))\leq n-1-\varepsilon,
\end{equation}
for some $\varepsilon>0$ depending only on $\alpha, \beta$.
In the more general case of integral functionals of the type \eqref{Mod2}, the theory of regularity is much less developed. The first regularity result established in the broader context of integral energies of the type $\eqref{Mod2}$ was accomplished by F. H. Lin and R. V. Kohn in 1999, indeed  in \cite{LK} they proved, for minimal configurations $(E,u)$ of \eqref{Mod2}, that $\mathcal {H}^{n-1}(\Sigma(E))=0$. The assumptions made by Lin and Kohn to achieve such regularity results require $C^k$ differentiability for $F$, $G$ and $\Psi$ as they appear in \eqref{Mod2}, where $k\geq 2$ and $F$, $G$ grow quadratically with respect to the gradient variable.  On the other hand,
nothing is proved concerning the Hausdorff dimension of $\Sigma(E)$ like in \eqref{redu}.\\
\indent We also point out that the problem can be set in a non-quadratic framework as well. This instance is less studied and only few regularity results are available (see \cite{AR, CEL, CFP, CFP2, E, FF, Lam}).\\
\indent In some recent papers, such as \cite{EL1,EL2}, the Hausdorff dimension estimate of $\Sigma(E)$ has been attained, significantly relaxing the differentiability assumption required on $F$ and $G$. Indeed, in \cite{EL2}, only H\"older continuous dependence of $F$ and $G$ with respect to $x$ and $u$ is necessary. However, it is worth noting that the aforementioned result is demonstrated under the assumption that $\Psi(\nu)=|\nu|$, representing the conventional perimeter.\\
\indent In this paper we will deal with the anisotropic case. 
Anisotropic surface energies manifest in various physical phenomena, such as crystal formation (refer to \cite{BNP1,BNP2}), liquid droplets (see \cite{CNT,DHV,FMA,LH,MV,Tay}), and capillary surfaces (see \cite{DFM, DM}). F. J. Almgren was a pioneer in investigating the regularity of surfaces that minimize anisotropic variational problems in his seminal paper \cite{A}. \
Early studies in this field were primarily conducted within the framework of varifolds and currents. While these results can be applied to surfaces of any codimension, they necessitate relatively stringent regularity assumptions on the integrands of the anisotropic energies, as outlined in \cite{Bom,SS}.\
More recently, the regularity assumptions on the integrands $\Psi$ of the anisotropic energies have been relaxed, as highlighted in \cite{DS,F}, where it is assumed that $\Psi(x,\cdot)$ is of class $C^1$ and $\Psi(\cdot,\xi)$ is H\"older continuous.\\
\indent
In this context, it is worth mentioning a very recent paper, \cite{Simm}, which establishes the regularity result for quasi-minimizers of anisotropic surface energies within the class of sets of finite perimeter, under the assumption of H\"older continuous dependence of $\Psi$ on $x$. This outcome is derived within the scope of ellipsoidal variational energies, as detailed in \eqref{Asurface}. Notably, surface energies of this specific form were initially introduced in a paper \cite{T} by J. Taylor.
In more detail, in the case that the elliptic integrand is given by $\Psi(x,\nu)=\langle A(x)\nu,\nu\rangle^{1/2}$, where $A(x)=\bigl(a_{ij}(x)\bigr)^n_{i,j=1}$ is an elliptic and bounded matrix, the surface energy  takes the form
\begin{equation}\label{Asurface}
\mathbf{\Phi}_A(E;G):=\int_{G\cap\partial^*E}\langle A(x)\nu_E,\nu_E\rangle^{1/2} \;d\mathcal{H}^{n-1}(x).
\end{equation}
We assume that $A$ is uniformly elliptic, that is there exist two constants $0<\lambda\leq \Lambda<+\infty$ such that
\begin{equation}\label{uell0}
\lambda|\xi|^2\leq\langle A(x)\xi,\xi\rangle\leq \Lambda |\xi|^2,\quad\forall x\in\Omega,\,\forall\xi\in\R^n.
\end{equation}
We require that $A$ is H\"older continuous with exponent $\mu\in (0,1]$, that is
\begin{equation}\label{Hol0}
[A]_{C^{\mu}(\Omega)}=\sup_{\substack{x\neq y\\x,y\in \Omega}}\frac{|A(x)-A(y)|}{|x-y|^{\mu}}<+\infty.
\end{equation}
To avoid excessive technicalities, we assume that the bulk energy follows a Dirichlet-type distribution, although the outcome could be generalized to functionals of type \eqref{Mod2}. Given a bounded open set $\Omega \subset \R^n$, we consider the following functional:
\begin{align}\label{DefFA}
&\mathcal{F}_A(E,u;\Omega)=\int_{\Omega}\sigma_E|\D u|^2\,dx+\mathbf{\Phi}_A(E;\Omega),
\end{align}
where $\sigma_E=\alpha \mathbbm{1}_E+\beta \mathbbm{1}_{\Omega\setminus E}$, $0<\alpha<\beta$ and $E\subset\Omega$. 
The achieved regularity result is presented within the scope of local almost-minimizers. This makes it applicable in a variety of concrete applications, as we will demonstrate, for example, in the case of constrained problems.  The following definition naturally arises in several problems from material sciences (see for example \cite{AB,EF,Lin,LK,Ma}, compare also with  \cite[Definition 2.2]{Simm}).
\begin{Def}[$(\kappa,\mu)$-minimizers]
\label{k-mu}
Let $U\Subset\Omega$. The energy pair $(E,u)$ is a $(\kappa,\mu)$-minimizer in $U$ of the functional $\mathcal {F}_{A}$, defined in \eqref{DefFA}, if for every $B_r(x_0)\subset U$ 
\begin{equation*}
\mathcal{F}_A(E,u;B_r(x_0))\leq\mathcal{F}_A(F,v;B_r(x_0))+\kappa |E\vartriangle F|^{\frac{n-1}{n}+\frac{\mu}{n}}, 
\end{equation*}
whenever $(F,v)$ is an admissible test pair, namely, $F$ is a set of finite perimeter with $F\vartriangle E\Subset B_r(x_0)$ and $v-u\in H^1_0(B_r(x_0))$.
\end{Def}
The main theorem proved in this paper is the following.
\begin{Thm}
\label{Teorema principale}
Let $(E,u)$ be a $(\kappa,\mu)$-minimizer of $\mathcal{F}_A$. Then
\begin{itemize}
\item[a)] there exists a relatively open set $\Gamma\subset \partial E$ such that $\Gamma$ is a $C^{1,\sigma}$-hypersurface, for all $0<\sigma<\frac{\mu}{2}$,
\item[b)] there exists $\varepsilon >0$ depending on $n,\alpha,\beta$ such that $$\mathcal{H}^{n-1-\varepsilon}((\partial E\setminus \Gamma)\cap \Omega)=0.$$
\end{itemize}
\end{Thm}
We outline the strategy adopted to prove this result. In the regularity theory for $\Lambda$-minimizers of the perimeter, the regular part $\Gamma$ of the boundary of $E$ is detected by the points that have a uniformly small excess in some ball (see Definition \ref{ExcessDef}). A decay relation for the excess plays a crucial role, as it triggers an iteration argument that shows that the unitary normal vector varies continuously along $\Gamma$, thus ensuring its smoothness.\\
\indent For our problem, it is not possible to prove a decay relation for the excess without considering the interaction between the surface and the bulk energy. Indeed, as outlined in Section \ref{Proof of the main theorem}, if the excess of a point $x_0$ in $\dd E$ in some ball $B_r(x_0)$ is small, we are only able to prove an improvement relation for the excess, which involves the rescaled Dirichlet integral of $u$ in $B_r(x_0)$ as well.\\
\indent In this context, $\Gamma$ is defined as the collection of the points of the boundary of $E$ that are centers of balls $B_r(x_0)$ where the excess $\mathbf{e}(E,x_0,r)$  and the rescaled Dirichlet integral $\mathcal{D}_u(x_0,r)$ (see    \eqref{Diric-rescaled}) of $u$ are sufficiently small. In order to prove the smoothness of $\Gamma$, a decay relation for the sum of these quantities is required.\\
\indent A decay relation for the rescaled Dirichlet integral of $u$ around points of small excess is proved separately in Proposition \ref{DecayDirichlet}.\\
\indent A much finer argument is needed to establish an improvement relation for the excess (see Theorem \ref{MiglEccesso}). 
One of the key concepts enabling us to adapt the standard excess-decay arguments, commonly used in the context of perimeter minimizers, to the anisotropic setting is a specific change of variable $T_{x_0}$. This affine transformation, already used in \cite{DEST,JPS,Simm}, maps Wulff shapes of $\Phi_A$, which are ellipsoids, into balls $B_r(x_0)$ (see Section \ref{Notation and preliminaries}). 
We first prove a version of the excess improvement theorem for transformed couples $(\tilde{E},\tilde{u})=\big(T_{x_0}(E),u\circ T_{x_0}^{-1}\big)$, which are $\big(\kappa\lambda^{-\frac{n}{2}},\mu\big)$-minimizers of $\mathcal{F}_{x_0,A_{x_0}}$ (see \eqref{Fx0D}). The proof of the latter is carried out by contradiction and is based on a blow-up argument. In this step, we can benefit of using the classical perimeter instead of the anisotropic one around $x_0$, being $A_{x_0}(x_0)=I$ (see Proposition \ref{MimimalitàCasoParticolare}). The main ideas can be summarized as follows:
\begin{enumerate}
     \item Density lower and upper bounds on the perimeter (see Theorem \ref{Energy upper bound} and Theorem \ref{Density lower bound}) guarantee that around points $x_0$ of $\dd \tilde{E}$ with small excess, the boundary of $\tilde{E}$ almost coincides with the graph of a Lipschitz function $f$ (see Theorem \ref{LipApp}). Therefore, it is possible to apply the area formula directly along $\dd \tilde{E}$ up to a small error. The portion of the boundary that does not match is controlled by the excess at that scale.
     \item The function $f$ is quasi-harmonic. We need a quantitative estimate of its quasi-harmonicity by a power of the excess with an exponent greater than $\frac{1}{2}$. In this step, the minimality of the optimal couple $(\tilde{E},\tilde{u})$ and the first variation formulae play a crucial role. 
     \item The direction of improvement of the excess is detected by the unitary normal vector to the graph of $f$. By means of a reverse Poincaré inequality (see Theorem \ref{RevPI}), the excess at a smaller scale at $x_0$ is controlled by the flatness of $\dd \tilde{E}$ around $x_0$, which is in turn estimated by the excess via the good properties of $f$.
     \end{enumerate}
\indent The paper is divided in sections, which reflect the proof strategy. Section \ref{Notation and preliminaries} collects notation and
preliminary definitions. In Section \ref{Scaling and change of variables}, some invariance properties of the excess and minimality under the transformation $T_{x_0}$ and rescaling are proved. In Section \ref{Energy density estimates}, we establish density lower and upper bounds for the perimeter of $E$ and their consequences, which are the decay of the rescaled Dirichlet energy and the Lipschitz approximation theorem. Section \ref{Compactness for sequences of minimizers} is devoted to prove a compactness result for sequences of $(\kappa,\alpha)$-minimizers, which serves as a crucial tool for estimating the size of the singular set of $E$, as stated in Theorem \ref{Teorema principale}. Section \ref{Reverse Poincaré inequality} includes the reverse Poincaré inequality, which is the counterpart of the well-known Caccioppoli's inequality for weak solutions of elliptic equations. Section \ref{Proof of the main theorem} contains the proof of Theorem \ref{Teorema principale}. The main ingredients to achieve such a result are a first variation formula for the bulk energy of the functional $\mathcal{F}_{x_0,A_{x_0}}$ and two versions of the excess improvement theorem. Finally, Section \ref{From costrained to penalized problem} deals with the application of the regularity result to a volume constrained problem via a penalization argument.

\section{Notation and preliminaries}
\label{Notation and preliminaries}
In the rest of the paper we will write $\langle \xi, \eta \rangle$ for the inner product of vectors $\xi, \eta \in \mathbb{R}^n$, and consequently $|\xi|:=\langle \xi, \xi \rangle^{\frac 12}$ will be the corresponding Euclidean norm. As usual $\omega_n$ stands for the Lebesgue measure of the unit ball in $\R^n$. We denote by $\mathbbm{S}^{n-1}$ the unit sphere of $\R^n$.\\
\indent We will write $x=(x',x_n)$ for all $x\in \R^n$, where $x'\in\R^{n-1}$ collects the first $n-1$ components of $x$ and $x_n\in\R$ is its $n$-th component.
Accordingly, we denote by $\D'=(\partial_{x_1},\dots,\partial_{x_{n-1}})$ the gradient with respect to the first $n-1$ components.\\  
\indent The $n$-dimensional ball in $\R^n$ with center $x_0$ and radius $r>0$ is denoted as
$$B_r(x_0)=\{x\in\R^{n}: |x-x_0|<r\}.$$ 
If $x_0=0$, we simply write $B_r$ in place of $B_r(x_0)$. The $(n-1)$-dimensional ball in $\R^{n-1}$ with center $x'_0$ and radius $r>0$ is denoted by 
\begin{equation*}
    \mathbf{D}_r(x'_0)=\{x'\in\R^{n-1}: |x'-x'_0|<R\}.
\end{equation*}
If $u$ is integrable in $B_R(x_0)$, we set
\begin{equation*}
    u_{x_0,r}=\frac{1}{\omega_n r^n}\int_{B_r(x_0)} u\,dx = \dashint_{B_{r}(x_0)} u\,dx.
\end{equation*}
If $E\subset \R^n$ and $t\in [0,1]$, the set of points of $E$ of density $t$ is defined as
\begin{equation}
E^{(t)}=\left\{x\in \R^n:\ |E\cap B_r(x)|=t|B_r(x)|+o(r^n)\text{ as }r\rightarrow 0^+\right\}.
\end{equation}
Given a Lebesgue measurable set $E\subset \R^n$, we say that $E$ is of locally finite perimeter if there exists a $\R^n$-valued Radon measure $\mu_E$ (called the Gauss-Green measure of $E$) such that
\begin{equation*}
\int_{E}\nabla \phi\ dx=\int_{\R^n}\phi \, d\mu_E,\quad\forall \phi \in C^1_c(\R^n).
\end{equation*}
Moreover, we denote the perimeter of $E$ relative to $G\subset \R^n$ by $P(E;G)=|\mu_E|(G)$.\\
 The support of $\mu_E$ can be characterized by
\begin{equation}\label{support}
\text{spt}\mu_E=\bigl\{x\in \R^n: 0<|E\cap B_r(x)|<\omega_n r^n, \,\forall r>0\bigr\},
\end{equation}
(see \cite[Proposition 12.19]{Ma}). It holds that $\text{spt}\mu_E\subset \partial E$.
 If $E$ is of locally finite perimeter, then the {\it reduced boundary} $\partial^*E$ of $E$ is the set of those $x\in\R^n$ such that
\begin{equation}\label{RB}
\nu_E(x):=\lim_{r\rightarrow 0^+}\frac{\mu_E(B_r(x))}{|\mu_E|(B_r(x))}
\end{equation}
exists and belongs to $\mathbb S^{n-1}$. In the following, the
topological boundary $\partial E$ must be understood by considering the representative $E^{(1)}$ of $E$, for which it holds that $\overline{\partial^*E}=\partial E$.

The properties of the matrix $A$ in the definition of the anisotropic perimeter \eqref{Asurface} guarantee that $ \mathbf{\Phi}_A(E;F)$ is comparable to the classical perimeter, as observed in \cite{Simm}. 

\begin{Rem}[Comparability to perimeter] 
$\mathbf{\Phi}_A(E;\cdot)$ is comparable to $P(E;\cdot)$, since for Borel sets $F\subset\R^n$, by the uniform ellipticity of $A$, it follows that
\begin{equation}\label{comparability}
\lambda^{1/2}P(E;F)\leq \mathbf{\Phi}_A(E;F)\leq \Lambda^{1/2}P(E;F).
\end{equation}
If $A$ equals the identity matrix $I$, we have the isotropic case $\mathbf{\Phi}_A(E;\cdot)=P(E;\cdot)$.
\end{Rem}

It will be useful in the sequel to build comparison sets by replacing regions within an open set. The anisotropic perimeter can be split as in the isotropic case.
\begin{Prop}[Comparison by replacement]\label{compreplace}
If $E$ and $F$ are sets of locally finite perimeter in $\R^n$ and $G$ is an open set of finite perimeter in $\R^n$ such that
\begin{equation}\label{nullinter}
\mathcal{H}^{n-1}(\partial^* G\cap \partial^*E)=\mathcal{H}^{n-1}(\partial^* G\cap \partial^*F)=0,
\end{equation}
then the set defined by
\begin{equation}\label{replacement}
F_0=(F\cap G)\cup (E\setminus G)
\end{equation}
is a set of locally finite perimeter in $\R^n$. Moreover if $G \Subset U$ and $U$ is an open subset of $\R^n$, then
\begin{equation}
\mathbf{\Phi}_A(F_0;U)=\mathbf{\Phi}_A (F;G)+\mathbf{\Phi}_A(E;U\setminus \overline{G})+\mathbf{\Phi}_A(G;E^{(1)}\vartriangle F^{(1)})
\end{equation}
\end{Prop}
\begin{proof}
The proof can be easly obtained from \cite[Theorem 16.16]{Ma}. Its details can be found in \cite[Proposition 4.3]{Simm}.
\end{proof}

In the following, for $R>0$ and  $\nu\in\mathbbm{S}^{n-1}$, we will denote the cylinder centered in $x_0$ with radius $R$ oriented in the direction $\nu$ by
\begin{equation*}
\mathbf{C}_{R}(x_0,\nu):=x_0+\{y\in\R^n\,:\,|\langle y,\nu \rangle|<R,\,|y-\langle y,\nu\rangle \nu|<R\},
\end{equation*}
and the cylinder of radius $R$ oriented in the direction $e_n$ with height 2 by
\begin{equation*}
\mathbf{K}_{R}(x_0):=
\mathbf{D}_{R}(x_0')\times (-1,1).
\end{equation*}
In the following, for simplicity of notation we will write $\mathbf{C}_R=\mathbf{C}_R(0,e_n)$
and $\mathbf{K}_{R}=\mathbf{K}_{R}(0)$.\\
\indent In addition we introduce some usual quantities involved in regularity theory.
\begin{Def}[Excess]\label{ExcessDef}
Let $E$ be a set of locally finite perimeter, $x\in\dd E$, $r>0$ and $\nu\in \mathbbm S^{n-1}$. We define:
\begin{itemize}
\item
the \textbf{cylindrical excess} of $E$ at the point $x$, at the scale $r$ and with respect to the direction $\nu$, as
\begin{equation*}
\ecc(E,x,r,\nu):=\frac{1}{r^{n-1}}\int_{\mathbf{C}_r(x,\nu)\cap\dd^*E}\frac{|\nu_E-\nu|^2}{2}d\,\mathcal{H}^{n-1}=\frac{1}{r^{n-1}}\int_{\mathbf{C}_r(x,\nu)\cap\dd^*E}[1-\langle\nu_E,\nu\rangle]\,d\mathcal{H}^{n-1}.
\end{equation*}
\item the \textbf{spherical excess} of $E$ at the point $x$, at the scale $r$ and with respect to the direction $\nu$, as
\begin{equation*}
{\mathbf e}(E,x,r,\nu):= \frac{1}{r^{n-1}}\int_{\partial^* E\cap B_r(x)}\frac{|\nu_E-\nu|^2}{2}d\mathcal H^{n-1}.
\end{equation*}
\item the \textbf{spherical excess} of $E$ at the point $x$ and at the scale $r$, as
\begin{equation*}
{\mathbf e}(E,x,r):=\min_{\nu \in \mathbbm S^{n-1}}{\mathbf e}(E,x,r,\nu).
\end{equation*}
\end{itemize} 
We omit the dependence on the set when it is clear from the context.
\end{Def}

\section{Scaling and change of variables}

Given a symmetric positive matrix $A$, the $D$ matrix of his eigenvalues and the matrix $V$  of ortonormal eigenvectors, we have $A=VDV^{-1}$. Accordingly we define $A^{1/2}=VD^{1/2}V^{-1}$. Being $A^{-1/2}A^{1/2}=I$, the anisotropic perimeter of $E$ coincides with the standard perimeter of the image of $E$ under the affine change of variable $y=A^{-1/2}x$ up to the scaling factor $\mathrm{det}\big(A^{-1/2}\big)$. Localizing this argument freezing the matrix $A^{-1/2}(x)$ in a point $x_0\in\dd E$, we define the affine change of variables
\begin{equation}\label{affinechange}
    T_{x_0}(x)=A^{-1/2}(x_0)(x-x_0)+x_0, \quad \quad T^{-1}_{x_0}(y)=A^{1/2}(x_0)(y-x_0)+x_0, \quad\forall x,y\in\R^n,
\end{equation}
and the matrix-valued function
\begin{equation*}
    A_{x_0}(y):=A^{-\frac{1}{2}}(x_0)A(T^{-1}_{x_0}(y))A^{-\frac{1}{2}}(x_0), \quad\forall y\in\R^n,
\end{equation*}
which satisfies $A_{x_0}(x_0)=I$. 
It can be easily verified that the set $T^{-1}_{x_0}(B_r(x_0))$,
is the Wulff shape of $\mathbf{\Phi}_{A(x_0)}$. Moreover the following inclusions hold:
\begin{equation*}
    B_{\lambda^{\frac{1}{2}} r}(x_0)\subset T^{-1}_{x_0}(B_r(x_0))\subset B_{\Lambda^{\frac{1}{2}} r}(x_0),
\end{equation*}
for any $r>0$ and $x_0\in\R^n$. Under the affine change of variable $T_{x_0}$, the minimality with respect to the functional $\mathcal{F}_A$ will be rephrased through the following functional
\begin{equation}
\label{Fx0D}
\mathcal{F}_{x_0,D}(E,u;\Omega)=\int_{\Omega}\sigma_E |\D uA^{-\frac{1}{2}}(x_0)|^2\,dx+\mathbf{\Phi}_{D}(E;\Omega),
\end{equation}
(see Proposition \ref{Invariance}).

In the sequel, we collect two invariance properties of $(\kappa,\mu)$-minimizers under the transformation $T_{x_0}$ and rescaling.
\label{Scaling and change of variables}
\begin{Prop}[Invariance of almost-minimizers under $T_{x_0}$]
\label{Invariance}
Let $(E,u)$ be a $(\kappa, {\mu})$-minimizer of $\mathcal{F}_A$ in $\Omega$ and $x_0\in\Omega$. Then $(T_{x_0}(E),u\circ T_{x_0}^{-1})$ is a $(\kappa\lambda^{-\frac{n}{2}}, {\mu})$-minimizer of $\mathcal{F}_{x_0,A_{x_0}}$ in the balls $B_r(z)$ such that $B_{({\Lambda}/{\lambda})^{\frac12}\;r}(z)\subset T_{x_0}(\Omega)$.
\end{Prop}

\begin{proof}
We use the notation $E_{0}:=T_{x_0}(E)$ and $u_{0}:=u\circ T_{x_0}^{-1}$. Let $B_{r}(z)\subset T_{x_0}(\Omega)$ be such that $B_{({\Lambda}/{\lambda})^{\frac12}\;r}(z)\subset T_{x_0}(\Omega)$ and $(F_{0},v_{0})$ be an admissible test pair, i.e. $F_{0}$ is a set of finite perimeter with $F_{0}\vartriangle E_{0} \Subset B_{r}(z)$ and $v_{0}-u_{0}\in H^1_0(B_r(z))$. First we notice that, setting $F=T^{-1}_{x_{0}}(F_{0})$, as in \cite[Proposition 4.1]{Simm}, we have
\[
E \vartriangle F \Subset T^{-1}_{x_{0}}(B_{r}(z)) \subset B_{\Lambda^{\frac{1}{2}}r}\big(T^{-1}_{x_{0}}(z)\big)\subset \Omega,
\]
where the last condition is satisfied because $r<\Lambda^{-\frac{1}{2}}\mathrm{dist}\big(T^{-1}_{x_0}(z),\dd\Omega\big)$. 
Moreover, for $v=v_{0}\circ T_{x_{0}}$, we have $v-u\in H_{0}^{1}\big(T^{-1}_{x_{0}}(B_{r}(z))\big)$ whose extension to zero in $B_{\Lambda^{\frac{1}{2}}r}\big(T^{-1}_{x_{0}}(z)\big)$, denoted again by $v-u$, belongs to $H_{0}^{1}\big(B_{\Lambda^{\frac{1}{2}}r}\big(T^{-1}_{x_{0}}(z)\big)\big)$.
It follows, by the hypothesis of $(\kappa, {\mu})$-minimality of the pair $(E,u)$, that
\begin{equation}
\label{eq:01-minimality}
\mathcal{F}_A\left(E,u;B_{\Lambda^{\frac{1}{2}}r}\big(T^{-1}_{x_{0}}(z)\big)\right)\leq\mathcal{F}_A\left(F,v;B_{\Lambda^{\frac{1}{2}}r}\big(T^{-1}_{x_{0}}(z)\big)\right)+\kappa |E\vartriangle F|^{\frac{n-1+\mu}{n}}.
\end{equation}
This simplifies to
\begin{equation}
\label{eq:02-minimality}
\mathcal{F}_A\left(E,u;T^{-1}_{x_{0}}(B_{r}(z)))\right)\leq\mathcal{F}_A\left(F,v;T^{-1}_{x_{0}}(B_{r}(z))\right)+\kappa |E\vartriangle F|^{\frac{n-1+\mu}{n}}.
\end{equation}
We now calculate, using \cite[formula (4.9)]{Simm} and the change of variables $y=T_{x_{0}}(x)$,
\[
\begin{split}
\mathcal{F}_{x_0,A_{x_0}}(E_{0},u_{0};B_{r}(z)) = & \int_{B_{r}(z)}\sigma_{E_{0}} |\D u_{0}A^{-\frac{1}{2}}(x_0)|^2\,dy+\mathbf{\Phi}_{A_{x_0}}(E_{0};B_{r}(z)) \\
& = \mathrm{det}(A^{-\frac{1}{2}}(x_0)) \left(  \int_{T_{x_{0}}^{-1}(B_{r}(z))}\sigma_{E} |\D u|^2\,dx + \mathbf{\Phi}_{A}\left(E;T^{-1}_{x_{0}}(B_{r}(z))\right)  \right) \\
& = \mathrm{det}(A^{-\frac{1}{2}}(x_0))  \mathcal{F}_A\left(E,u;T^{-1}_{x_{0}}(B_{r}(z)))\right) \\
& \leq  \mathrm{det}(A^{-\frac{1}{2}}(x_0)) \left( \mathcal{F}_A\left(F,v;T^{-1}_{x_{0}}(B_{r}(z))\right)+\kappa |E\vartriangle F|^{\frac{n-1+\mu}{n}}\right) ,
\end{split}
\]
where with a slight abuse of notation we have denoted $\sigma_{E_{0}}={\alpha}\mathbbm{1}_{E_{0}}+ {\beta}\mathbbm{1}_{T_{x_{0}}(\Omega)\setminus E_{0}}$.
The result follows by observing that $\mathrm{det}(A^{-\frac{1}{2}}(x_0)) \leq \lambda^{-\frac{n}{2}}$ and 
\[
 \mathrm{det}(A^{-\frac{1}{2}}(x_0))  \mathcal{F}_A\left(F,v;T^{-1}_{x_{0}}(B_{r}(z))\right) = \mathcal{F}_{x_0,A_{x_0}}(F_{0},v_{0};B_{r}(z)).  
\]
\end{proof}

\begin{Prop}[Scaling of $(\kappa, {\mu})$-minimizers]

\label{PropScaling}
For $x_0\in\Omega$ and $r>0$, let $(E,u)$ be a 
$(\kappa,\mu)$-minimizer of $\mathcal{F}_{x_0,A_{x_0}}$ in $\Omega$  {(or a 
$(\kappa, \mu)$-minimizer of $\mathcal{F}_{A}$)}. 
Then $\big(\Psi_{x_0,r}(E),r^{-\frac{1}{2}}u\circ \Psi_{x_0,r}^{-1}\big)$ is a $(\kappa r^{\mu}, \mu)$-minimizer of $\mathcal{F}_{x_0,A_{x_0}\circ\Psi_{x_0,r}^{-1}}$ in $\Psi_{x_0,r}(\Omega)$, 
 {(or respectively a $(\kappa r^{\mu},\mu)$-minimizer of $\mathcal{F}_{A}$ in $\Psi_{x_0,r}(\Omega)$)}, where
\begin{equation*}
\Psi_{x_0,r}(x):=\frac{x-x_0}{r}, \quad\forall x\in\R^n.
    \end{equation*}
\end{Prop}

\begin{proof}
    Let $B_s(z)\subset \Psi_{x_0,r}(\Omega)$. Applying the change of variables $y=\Psi_{x_0,r}(x)$, we deduce that
    \begin{align*}
        & \mathcal{F}_{x_0,A_{x_0}\circ\Psi_{x_0,r}^{-1}}\big(\Psi_{x_0,r}(E),u\circ \Psi_{x_0,r}^{-1};B_s(z)\big)\\
        & =\int_{B_s(z)}\sigma_{\Psi_{x_0,r}(E)}|\D (u\circ \Psi_{x_0,r}^{-1})A^{-\frac{1}{2}}(x_0)|^2\,dy+\mathbf{\Phi}_{A_{x_0}\circ\Psi_{x_0,r}^{-1}}(\Psi_{x_0,r}(E);B_s(z))\\
        & =\frac{1}{r^{n-1}}\bigg(\int_{B_{rs}(x_0+rz)}\sigma_{E}|\D u A^{-\frac{1}{2}}(x_0)|^2\,dx+\mathbf{\Phi}_{A_{x_0}}(E;B_{rs}(x_0+rz))\bigg)\\
        & =\frac{1}{r^{n-1}}\mathcal{F}_{x_0,A_{x_0}}(E,u;B_{rs}(x_0+rz)).
    \end{align*}
    Let $(F,v)$ be such that $F$ is a set of finite perimeter, $F\vartriangle \Psi_{x_0,r}(E)\Subset B_s(z)$ and $v\in (u\circ \Psi^{-1}_{x_0,r})+H^1_0(B_s(z))$. It holds that $ \Psi_{x_0,r}^{-1}(F)\vartriangle E\Subset B_{rs}(x_0+rz)\subset\Omega$ and $v\circ\Psi_{x_0,r}-u\in H^1_0(B_{rs}(x_0+rz))$. Using the $(\kappa,\mu)$-minimality of $(E,u)$ we get
    \begin{align*}
        & \mathcal{F}_{x_0,A_{x_0}\circ\Psi_{x_0,r}^{-1}}\big(\Psi_{x_0,r}(E),u\circ \Psi_{x_0,r}^{-1};B_s(z)\big)=\frac{1}{r^{n-1}}\mathcal{F}_{x_0,A_{x_0}}(E,u;B_{rs}(x_0+rz))\\
        & \leq \frac{1}{r^{n-1}}\big[\mathcal{F}_{x_0,A_{x_0}}\big(\Psi_{x_0,r}^{-1}(F),v\circ\Psi_{x_0,r};B_{rs}(x_0+rz)\big)+\kappa|\Psi_{x_0,r}(E)\vartriangle\Psi_{x_0,r}(F)|^{\frac{n-1+\mu}{n}}\big]\\
        & =\mathcal{F}_{x_0,A_{x_0}\circ\Psi_{x_0,r}^{-1}}(F,v;B_s(z))+{\kappa r^{\mu}|E\vartriangle F|^{\frac{n-1+\mu}{n}}},
    \end{align*}
    which means that $\big(\Psi_{x_0,r}(E),r^{-\frac{1}{2}}u\circ \Psi_{x_0,r}^{-1}\big)$ is a $(\kappa r^{\mu}, \mu)$-minimizer of $\mathcal{F}_{x_0,A_{x_0}\circ\Psi_{x_0,r}^{-1}}$ in $\Psi_{x_0,r}(\Omega)$.
\end{proof}

The following proposition will be useful in the proof of Theorem \ref{MiglEccesso}, where we need to compare the excess of the transformed couple and of the original couple.

\begin{Prop}[Comparability of the excess under change of variable $T_{x_0}$ and radius]\label{CompExcess}
There exists a positive constant $c_1=c_1(n,\lambda,\Lambda)$ such that if $E$ is a set of locally finite perimeter and $x_0\in\dd E$, then, for any $r>0$,
\begin{equation*}
    c_1^{-1}\mathbf{e}\big(T_{x_0}(E),x_0,\Lambda^{-\frac{1}{2}}r\big)\leq \mathbf{e}(E,x_0,r)\leq c_1\mathbf{e}\big(T_{x_0}(E),x_0,\lambda^{-\frac{1}{2}}r\big).
\end{equation*}
    \begin{proof}
        For $r>0$ and $\nu\in\mathbbm{S}^{n-1}$, we define the ellipsoidal excess at the point $x_0$, at the scale $r$ and with respect to the direction $\nu$ the following quantity:
        \begin{equation*}
{\mathbf e}_W(E,x_0,r,\nu):=\min_{\nu \in \mathbbm S^{n-1}}\frac{1}{r^{n-1}}\int_{T^{-1}_{x_0}(B_r(x_0))\cap\dd^*E}\frac{|\nu_E-\nu|^2}{2}d\mathcal H^{n-1}.
\end{equation*}
Leveraging the inclusions 
\begin{equation*}
T^{-1}_{x_0}\big(B_{\Lambda^{-\frac{1}{2}}r}(x_0)\big)\subset B_r(x_0)\subset T^{-1}_{x_0}\big(B_{\lambda^{-\frac{1}{2}}r}(x_0)\big),
\end{equation*}
we infer that
\begin{equation}\label{equa1}
c^{-1}\mathbf{e}_W\big(E,x_0,\Lambda^{-\frac{1}{2}}r,\nu\big)\leq \mathbf{e}(E,x_0,r,\nu)\leq c\,\mathbf{e}_W\big(E,x_0,\lambda^{-\frac{1}{2}}r,\nu\big),
\end{equation}
for some positive constant $c=c(n,\lambda,\Lambda)$. 
In \cite[Proposition 5.1]{Simm}, it is shown that
\begin{equation*}
    c^{-1}\mathbf{e}\Bigg(T_{x_0}(E),x_0,s,\frac{A^{\frac{1}{2}}(x_0)\nu}{|A^{\frac{1}{2}}(x_0)\nu|}\Bigg)\leq \mathbf{e}_W(E,x_0,s,\nu)\leq c\,\mathbf{e}\Bigg(T_{x_0}(E),x_0,s,\frac{A^{\frac{1}{2}}(x_0)\nu}{|A^{\frac{1}{2}}(x_0)\nu|}\Bigg),
\end{equation*}
for any $s>0$. Using the previous inequalities from below for $s=\Lambda^{-\frac{1}{2}}r$ and from above for $s=\lambda^{-\frac{1}{2}}r$, and inserting them in \eqref{equa1}, we get
\begin{equation*}
    c^{-1}\mathbf{e}\Bigg(T_{x_0}(E),x_0,\Lambda^{-\frac{1}{2}}r,\frac{A^{\frac{1}{2}}(x_0)\nu}{|A^{\frac{1}{2}}(x_0)\nu|}\Bigg)\leq \mathbf{e}(E,x_0,r,\nu) \leq c\,\mathbf{e}\Bigg(T_{x_0}(E),x_0,\lambda^{-\frac{1}{2}}r,\frac{A^{\frac{1}{2}}(x_0)\nu}{|A^{\frac{1}{2}}(x_0)\nu|}\Bigg).
\end{equation*}
Minimizing over $\nu\in\mathbbm{S}^{n-1}$, we obtain the thesis.
    \end{proof}
\end{Prop}

\section{Energy density estimates}
\label{Energy density estimates}
The main goal of this section is to prove density lower and upper bounds for the perimeter of a $(\kappa,\mu)$-minimizer. As consequences, the decay of the associate rescaled Dirichlet energy, defined for $u\in H^1(B_r(x_0))$ as
\begin{equation}
\label{Diric-rescaled}
\mathcal{D}_u(x_0,r):=\frac{1}{r^{n-1}}\int_{B_r(x_0)}|\D u|^2\,dx,
\end{equation}
and the Lipshitz approximation theorem will follow.
In view of this aim, we mention a result stating a decay estimate for elastic minima around points where either the density of $E$ is close to $0$ or $1$, or the
set $E$ is asymptotically close to a hyperplane. We address the reader to \cite[Proposition 2.4]{FJ} for the proof.
\begin{Lem}
\label{Lemma decadimento 1}
Let $(E,u)$ be a {$(\kappa,\mu)$}-minimizer of the functional ${\mathcal F_A}$. There exists $\tau_0\in(0,1)$ such that the following statement is true: for all $\tau \in (0,\tau_0)$ there exists $\varepsilon_0=\varepsilon_0(\tau)>0$ such that if $B_r(x_0)\Subset \Omega$ and one of the following conditions holds:
\begin{itemize}
\item[\emph{(i)}]$ |E\cap B_r(x_0)|<\varepsilon_0 |B_r(x_0)|$,
\item[\emph{(ii)}]$ |B_r(x_0)\setminus E|<\varepsilon_0 |B_r(x_0)|$,
\item[\emph{(iii)}] There exists a halfspace $H$ such that $\frac{\left|(E\vartriangle H)\cap B_r(x_0)\right|}{|B_r(x_0)|}<\varepsilon_0$,
\end{itemize}
then
\begin{equation*}
\mathcal{D}_u(x_0,\tau r)\leq c_2\tau\mathcal{D}_u(x_0,r),
\end{equation*}
for some positive constant $c_2=c_2\big(n,\alpha,\beta\big)$.
\end{Lem}

The second result we want to mention, which will be used later, provides an upper bound for the whole energy $\mathcal{F}_A$ on balls. The proof is rather standard and we address the reader to \cite[Theorem 3]{EL2} for the details. Here we just give a sketch of the proof, underlining the only points where the presence of the anisotropy entails different computations.
\begin{Thm}[Energy upper bound]
\label{Energy upper bound}
Let $(E,u)$ be a $(\kappa,\mu)$-minimizer of $\mathcal{F}_A$ in $\Omega$. Then for every open set $U\Subset \Omega$ there exists a positive constant $c_3=c_3\big(n,\alpha,\beta,\Lambda,\kappa,\mu,U,\norm{\D u}_{L^2(\Omega)}\big)$ such that for every $B_r(x_0)\subset U$ it holds
\begin{equation*}
\mathcal{F}_A(E,u;B_r(x_0))\leq c_3 r^{n-1}.
\end{equation*}
\end{Thm}
\begin{proof}
Let $B_r(x_0)\subset U\Subset\Om$. Testing the minimality of $(E,u)$ with $(E\setminus B_r(x_0),u)$, we deduce that
\begin{equation}\label{mindup}
\mathcal{F}_A(E,u;\Omega)
\leq \mathcal{F}_A(E\setminus B_r(x_0),u;\Omega)+\kappa |E\vartriangle(E\setminus B_r(x_0))|^{\frac{\mu+n-1}{n}}
. 
\end{equation}
The only difference in our proof, compared to the isotropic case, is the use of the following formula concerning anisotropic perimeter and set operations. The latter follows from Proposition \ref{compreplace} applied with $F=\emptyset$, $G=B_r(x_0)$, that is
\begin{equation*}
\mathbf{\Phi}_A(E\setminus B_r(x_0);U)=\mathbf{\Phi}_A(E;U\setminus \overline{{B}_r(x_0)})+\mathbf{\Phi}_A(B_r(x_0);E).
\end{equation*}
Making $\mathcal{F}_A$ explicit and getting rid of the common terms in \eqref{mindup} we obtain  the following energy estimate on $B_r(x_0)\cap E$:
\begin{align}
\label{a36}
\int_{B_r(x_0)\cap E}(\beta -\alpha)|\D u|^2\,dx+\mathbf{\Phi}_A(E;B_r(x_0))
& \leq \mathbf{\Phi}_A(B_r(x_0);E)+\kappa 
|E\vartriangle(E\setminus B_r(x_0))|^{\frac{\mu+n-1}{n}}
\notag\\
& \leq \Lambda^{1/2}\mathcal{H}^{n-1}(\dd B_r(x_0))+c(n,\kappa)r^{n-1}\notag\\
&\leq c(n,\Lambda,\kappa)r^{n-1}.
\end{align}
Starting from this estimate the proof follows verbatim the argument used in \cite[Theorem 3]{EL2}, because henceforth only variations of the function $u$ are used and the perimeter is not involved anymore.
Indeed, using a blow up argument, it can be proved that there exist {$M=M(n,\alpha,\beta)>0$ and $\tau\in\big(0,\frac{1}{2}\big)$, depending on ${\Lambda}/{\lambda}$,} such that for every {$\delta\in(0,1)$} there exists $h_0\in\N$ such that, for any $B_r(x_0)\subset U$, we have
\begin{equation*}
\int_{B_r(x_0)}|\D u|^2\leq h_0r^{n-1} \quad\text{or}\quad \int_{B_{\tau r}(x_0)}|\D u|^2\,dx\leq M\tau^{n-\delta}\int_{B_r(x_0)}|\D u|^2\,dx,
\end{equation*}
from which the thesis follows.
\end{proof}
In the following lemma we show that the energy $\mathcal F_A$ decays ``fast'' in the balls where the perimeter of $E$ is ``small''. Lemma \ref{Lemma decadimento 1} is utilized in its proof, specifically in instances (i) and (ii).
\begin{Lem}
\label{Lemma decadimento 2} Let $(E,u)$ be a $(\kappa,\mu)$-minimizer in $\Omega$ of the functional ${\mathcal F_A}$.
 For every $\tau\in (0,1)$ there exists $\varepsilon_1=\varepsilon_1(\tau)>0$ such that, if $B_r(x_0)\subset \Omega$ and $P(E;B_r(x_0))<\varepsilon_1 r^{n-1}$, then
\begin{equation}\label{D}
\mathcal F_A(E,u;B_{\tau r}(x_0))\leq c_4 \bigl(\tau^n\mathcal F_A(E,u;B_r(x_0))+(\tau r)^{\mu+n-1}\bigr),
\end{equation}
for some positive constant $c_4=c_4\big(n,\alpha,\beta,\lambda,\Lambda,\kappa,\mu,\norm{\D u}_{L^2(\Omega)}\big)>0$ independent of $\tau$ and $r$.
\end{Lem}

\begin{proof}
Let $\tau\in(0,1)$ and $B_r(x_0)\subset\Omega$. Without loss of generality, we may assume that $\tau<\frac 12$. We rescale $(E,u)$ in $B_1$ by setting $E_r=\frac{E-x_0}{r}$ and $u_r(y)=r^{-\frac 12}{u(x_0+ry)}$, for $y\in B_1$. Applying Proposition \ref{PropScaling}, we have that $(E_r,u_r)$ is a $(\kappa r^{\mu},\mu)$-minimizer of $\mathcal F_{\tilde A}$ in 
$B_1$, where $\tilde{A}=A\circ\Psi_{x_0,r}^{-1}$. Observing that $r^{n-1}\mathcal F_{\tilde{A}}(E_r,u_r;B_{\tau })=\mathcal F_A(E,u;B_{\tau r }(x_0))$, we have to prove that there exists $\varepsilon_1=\varepsilon_1(\tau)$ such that, if $P(E;B_1)<\varepsilon_1$, then
\begin{equation}
\mathcal F_{\tilde A}(E,u;B_{\tau })\leq c_4 \big(\tau^n\mathcal F_{\tilde A}(E,u;B_1)+\tau^{\mu+n-1}r^\mu\big).
\end{equation}
For simplicity of notation we will still denote $E_r$ by $E$, $u_r$ by $u$ and $\tilde{A}$ by $A$. 
We note that, since $P(E;B_1)<\varepsilon_1$, by the relative isoperimetric inequality, either $|B_1\cap E|$ or $|B_1\setminus E|$ is small and thus Lemma \ref{Lemma decadimento 1} can be applied. We assume that $|B_1\setminus E|\leq |B_1\cap E|$, the other case being similar. 
By the coarea formula and the relative isoperimetric inequality we get
$$\int_\tau^{2\tau}\mathcal{H}^{n-1}(\partial B_\rho\setminus E)\,d\rho\leq|B_1\setminus E|\leq c(n) P(E;B_1)^{\frac{n}{n-1}}.$$
Therefore, we may choose $\rho\in (\tau,2\tau)$, independent of $n$, such that it holds $\mathcal{H}^{n-1}(\partial^*E\cap \partial B_{\rho})=0$ and
\begin{equation}
\label{eqq8}
\mathcal{H}^{n-1}(\partial B_\rho\setminus E)\leq \frac{c(n)}{\tau} P(E;B_1)^{\frac{n}{n-1}}\leq\frac{c(n)\varepsilon_1^{\frac{1}{n-1}}}{\tau} P(E;B_1).
\end{equation}
Now we test the minimality of $(E,u)$ with $(F_0,u)$, where $F_0:=E\cup B_{\rho}$.
We remark that, being $\mathcal{H}^{n-1}(\partial^*E\cap \partial B_{\rho})=0$, we can apply Proposition \ref{replacement} with $U=F=B_1$ and $G=B_\rho$, thus obtaining
\begin{equation}\label{eqq9}
 \mathbf{\Phi}_A(F_0;B_1)=\mathbf{\Phi}_A(E\cup B_{\rho};B_1)
    =\mathbf{\Phi}_A(E;B_1\setminus \overline{B_\rho})+\mathbf{\Phi}_A(B_{\rho};B_1\setminus E^{(1)}).
 \end{equation}
The $(\kappa r^{\mu},\mu)$-minimality of $(E,u)$ supplies
\begin{equation}
\label{eqq61}
\int_{B_1}\sigma_E|\D u|^2\,dx+\mathbf{\Phi}_A(E;B_1)
\leq \int_{B_1}\sigma_{F_0}|\D u|^2\,dx+\mathbf{\Phi}_A(F_0;B_1)+\kappa r^{\mu}{|E\vartriangle F_0|^{\frac{\mu+n-1}{n}}}.
\end{equation}
Using \eqref{eqq9} to get rid of the common perimeter terms and recalling that $E=E^{(1)}$, we deduce
\begin{equation}
\label{eqq7}
\int_{B_1}\sigma_E|\D u|^2\,dx+\mathbf{\Phi}_A(E;B_{\rho})
\leq \int_{B_1}\sigma_{F_0}|\D u|^2\,dx+\mathbf{\Phi}_A(B_{\rho};B_1\setminus E)+\kappa r^{\mu}|E\vartriangle F_0|^{\frac{\mu+n-1}{n}
}.
\end{equation}
Taking into account the comparability to the perimeter \eqref{comparability} and perimeter estimate \eqref{eqq8}, recalling that $\rho\in (\tau,2\tau)$ and getting rid of the common Dirichlet terms, we deduce:
\begin{align}
\label{eqq5}
\int_{B_{\tau}}\sigma_E|\D u|^2\,dx+\lambda^{1/2}P(E;B_{\tau})
& \leq \beta\int_{B_{2\tau}}|\D u|^2\,dx+\Lambda^{1/2}\mathcal{H}^{n-1}(\partial B_\rho\setminus E)+c(n,\kappa) r^{\mu}\tau^{\mu+n-1}\\
&\leq \beta\int_{B_{2\tau}}|\D u|^2\,dx+\frac{c(n)\Lambda^{1/2}}{\tau} \varepsilon_1^{\frac{1}{n-1}}P(E;B_1)+c(n,\kappa) r^{\mu}\tau^{\mu+n-1}.
\end{align}
Finally, we choose $\varepsilon_1$ such that

\begin{equation}
c(n)\Lambda^{1/2} \varepsilon_1^{\frac{1}{n-1}}\leq \tau^{n+1}\quad\mbox{and }\quad c(n)\varepsilon_1^{\frac{n}{n-1}} \leq \varepsilon_0(2\tau)|B_1|,
\end{equation}
where $\varepsilon_0$ is from Lemma \ref{Lemma decadimento 1}, thus getting
\begin{equation}
\int_{B_{2\tau}}|\D u|^2\,dx\leq 2^n c_2\tau^n\int_{B_1}|\D u|^2\,dx.
\end{equation}
From this estimates the result easily follows applying again the comparability to the perimeter.
\end{proof}
Taking advantage of the established results, we are able to deduce a density lower bound estimate for the perimeter of a $(\kappa,\mu)$-minimizer of $\mathcal{F}_A$.
\begin{Thm}[Density lower bound]
\label{Density lower bound}
Let $(E,u)$ be a $(\kappa,\mu)$-minimizer of $\mathcal{F}_A$ in $\Omega$ and $U\Subset \Omega$ be an open set. Then there exists a constant $c_5=c_5\big(n,\alpha,\beta,\lambda,\Lambda,\kappa,\mu,U,\norm{\D u}_{L^2(\Omega)}\big)>0$,
such that, for every $x_0\in \partial E$ and $B_r(x_0)\subset U$, it holds
\begin{equation}\label{DensityLB}
P(E;B_r(x_0))\geq c_5 r^{n-1}.
\end{equation}
Moreover, $\mathcal{H}^{n-1}((\partial E\setminus \partial^*E)\cap \Omega)=0$.
\end{Thm}
\begin{proof}
The proof matches that of \cite[Theorem 4]{EL2} exactly, given the comparability to the perimeter. We start by assuming that $x_0\in\displaystyle\dd^*E$. Without loss of generality we may also assume that $x_0=0$. Arguing by contradiction on \eqref{DensityLB}, by using Theorem \ref{Energy upper bound} and Lemma \ref{Lemma decadimento 2}, we can easily prove by induction (see \cite[Theorem 4]{EL2} for the details) that
\begin{equation}
\label{Relazione iterativa}
\mathcal{F}(E,u;B_{\sigma\tau^hr})\leq\varepsilon_1(\tau)\tau^{\mu h}(\sigma\tau^h r)^{n-1},
\end{equation}
where $\tau$ and $\sigma$ are sufficiently small and $\varepsilon_1$ is from Lemma \ref{Lemma decadimento 2}.
Starting from this, we deduce that
\begin{equation*}
\lim_{\rho\rightarrow 0^+}\frac{P(E;B_\rho)}{\rho^{n-1}}=\lim_{h\rightarrow +\infty}\frac{P(E; B_{\sigma\tau^hr})}{(\sigma\tau^hr)^{n-1}}\leq\lim_{h\rightarrow+\infty}2\varepsilon_1(\tau)\tau^{\mu h}=0,
\end{equation*}
which implies that $x_0\not\in\dd^*E$, that is a contradiction.
We recall that we chose the representative of $\dd E$ such that $\dd E=\overline{\dd\displaystyle^*E}$.
Thus, if $x_0\in\dd E$, there exists $(x_h)_{h\in\N}\subset\dd^*E$ such that $x_h\rightarrow x_0$ as $h\rightarrow+\infty$,
\begin{equation*}
P(E;B_r(x_h))\geq C r^{n-1}
\end{equation*}
and $B_r(x_h)\subset U$, for $h$ large enough. Passing to the limit as $h\rightarrow+\infty$, we get the result.
\end{proof}

\begin{Def}[Ahlfors regularity]\label{Ahlfors}
A Borel measure $\mathbb{\mu}$ on $\R^n$ is said to be $d$-Ahlfors regular if there exist two positive constants $c_A$ and $r_0$ such that
\begin{equation}\label{lowupperbound}
c_A^{-1}r^d\leq \mathbb{\mu}(B_r(x))\leq c_A r^d,
\end{equation}
for all $x\in \mathrm{spt}\mathbb{\mu}$ and $0<r<r_0$. According to the notation used in \cite{BEGTZ}, we denote
\begin{align}
\mathcal{A}(c_A,r_0):
=\Bigl\{& E\subset \R^n : E \text{ is a set of locally finite perimeter satisfying}
\\
&\partial E=\mathrm{spt}\mathbb{\mu}_E \text{ and its perimeter measure } 
|\mathbb{\mu}_E| \text{ is}
\\
&(n-1)\text{-Ahlfors regular with constants } r_0 \text{ and }c_A\Bigr\}.
\end{align}
\end{Def}
\begin{Rem}\label{DensityBounds} It is evident that Theorems \ref{Energy upper bound} and \ref{Density lower bound} ensure the belonging of the $(\kappa,\mu)$-minimizers of $\mathcal{F}_A$ to the class $\mathcal{A}(c_A,r_0)$, for some constant $c_A$ identified in such theorems. Naturally, for $x_0\in\Omega$ and $r>0$, the 
$\big(\kappa \lambda^{-\frac{n}{2}}r^\mu,\mu\big)$-minimizers of $\mathcal{F}_{x_0,A_{x_0}\circ\Psi_{x_0,r}^{-1}}$ obtained through the affine transformation $T_{x_0}$ and the scaling $\Psi_{x_0,r}^{-1}$ (see Proposition \ref{Invariance} and Proposition \ref{PropScaling}) belong to the class $\mathcal{A}\Big(c_A\big(\frac{\Lambda}{\lambda}\big)^{\frac{n}{2}},\frac{r_0}{r\Lambda^{\frac{1}{2}}}\Big)$. 
\end{Rem}

The next result of this section establishes that around the points of the boundary of the set where the excess is ``small'', the Dirichlet integral decays ``fast''. In its proof, Lemma \ref{Lemma decadimento 1} plays a crucial role in istance (iii).

\begin{Prop}[Decay of the rescaled Dirichlet integral]\label{DecayDirichlet}
For every $\tau\in(0,1)$ there exists $\varepsilon_2=\varepsilon_2(\tau)>0$ such that if $(E,u)$ is a $(\kappa,\mu)$-minimizer of $\mathcal{F}_A$ in $B_r(x_0)$, with $x_0\in\dd E$, and $\mathbf{e}(x_0,r)\leq\varepsilon_2$, 
then \begin{equation}
\mathcal{D}_u(x_0,\tau r)\leq c_6\tau\mathcal{D}_u(x_0,r),
\end{equation}
for some positive constant $c_6=c_6\big(n,\alpha,\beta,\norm{\D u}_{L^2(\Omega)}\big)$.
\end{Prop}
\begin{proof}
Applying a usual scaling argument, by Proposition \ref{PropScaling}, we assume by contradiction that for some $\tau\in(0,1)$ there exist two positive sequences $(\varepsilon_h)_{h\in\N}$ and $(r_h)_{h\in\N}$ and a sequence $((E_h,u_h))_{h\in\N}$ of $(\kappa r_h^\mu,\mu)$-minimizers of $\mathcal{F}_{A\circ\Psi_{x_0,r_h}^{-1}}$ in $B_1$ with equibounded energies such that $0\in\dd E_h$,
\begin{equation}
\label{Eqn 7}
{\mathbf e}(E_h,0,1)=\varepsilon_h\rightarrow 0 \quad\text{and}\quad \mathcal{D}_{u_h}(0,\tau)>\overline C\tau\mathcal{D}_{u_h}(0,1),
\end{equation}
for some positive constant $\overline C$ to be chosen.
Thanks to the energy upper bound (Theorem \ref{Energy upper bound}) and the compactness of $(E_h)_{h\in\N}$, we may assume that $E_h\rightarrow E$ in $L^1(B_1)$ and $0\in\dd E$. Since, by lower semicontinuity, the excess of $E$ at 0 is null, $E$ is a half-space in $B_1$, say $H$. In particular, for $h$ large, it holds
\begin{equation*}
|(E_h\vartriangle H)\cap B_1|<\varepsilon_0(\tau)|B_1|,
\end{equation*}
where $\varepsilon_0$ is from Lemma \ref{Lemma decadimento 1}, which gives a contradiction with the inequality \eqref{Eqn 7}, provided we choose $\overline{C}>c_2$, where $c_2$ is also from Lemma \ref{Lemma decadimento 1}.
\end{proof}

The last results also come as consequences of the density lower and upper bounds proved above. The height bound lemma is a standard step in the proof of regularity because it is one of the main ingredients to prove the Lipschitz approximation theorem. We remark that this is stated for $(\kappa r^\mu,\mu)$-minimizers of $\mathcal{F}_{A_{x_0}\circ\Psi_{x_0,r}}$, which are still Ahlfors regular (see Remark \ref{DensityBounds}). The proof of this result can be found in \cite[Theorem A.2]{BEGTZ}.

\begin{Lem}[Height bound]
\label{height bound}
For $x_0\in\Omega$ and $r>0$, let $(E,u)$ be a $(\kappa r^\mu,\mu)$-minimizer of $\mathcal{F}_{A_{x_0}\circ\Psi_{x_0,r}^{-1}}$ in $B_1$. There exist two positive constants $\varepsilon_3$ and $c_7$, depending on $n,\alpha,\beta,\lambda,\Lambda,\kappa,\mu,\norm{\D u}_{L^2(B_1)}$, such that if $0\in \partial E$ and
\begin{equation*}
{\mathbf e}(0,1,e_n)<\varepsilon_3,
\end{equation*}
then
\begin{equation*}
\sup_{y\in \partial E\cap B_{1/2}}|y_n-(x_0)_n|\leq c_7{\mathbf e}(0,1,e_n)^{\frac{1}{2(n-1)}}.
\end{equation*}
\end{Lem}

Proceeding as in \cite{Ma}, we state the following Lipschitz approximation lemma, which is a consequence of the height bound lemma. Its proof follows exactly as in \cite[Theorem A.3]{BEGTZ}. {It is a foundamental step in the long journey to the regularity because it provides a connection between the regularity theories for parametric and non-parametric variational problems. Indeed we are able to prove for $(\kappa r^\mu,\mu)$-minimizers that the smallness of the excess guaranties that $\partial E$ can be locally almost entirely covered by the graph of a Lipschitz function}.

\begin{Thm}[Lipschitz approximation]\label{LipApp}
For $x_0\in\Omega$ and $r>0$, let $(E,u)$ be a $(\kappa r^\mu,\mu)$-minimizer of $\mathcal{F}_{A_{x_0}\circ\Psi_{x_0,r}^{-1}}$ in $B_1$.
There exist two positive constants $\varepsilon_4$ and $c_8$, depending on\break $n,\alpha,\beta,\lambda,\Lambda,\norm{\D u}_{L^2(B_1)}$, such that if $0\in \partial E$ and
$$
{\mathbf e}(0,1,e_n)<\varepsilon _4,
$$
then there exists a Lipschitz function $f:\R^{n-1}\rightarrow \R$ such that
$$
\sup_{x'\in \R^{n-1}}|f(x')|\leq c_8{\mathbf e}(0,1,e_n)^{\frac{1}{2(n-1)}},\quad \norm{\nabla'f}_{L^{\infty}}\leq 1,
$$
and
$$
\mathcal{H}^{n-1}((\partial E \vartriangle \Gamma_f)\cap B_{1/2})\leq c_8 {\mathbf e}(0,1,e_n),
$$
where $\Gamma_f$ is the graph of $f$. Moreover,
$$
\int_{\mathbf{D}_\frac{1}{2}}|\nabla'f|^2\,dx'\leq c_8 {\mathbf e}(0,1,e_n).
$$
\end{Thm}

\section{Compactness for sequences of minimizers}
\label{Compactness for sequences of minimizers}
In this section we prove a standard compactness result for sequences of $(\kappa,\mu)$-minimizers. Given to positive constants $M_1$ and $M_2$, we set 
\begin{equation*}
B_{M_1,M_2}:=\bigl\{A\in C^{\gamma}(\R^n;\R^n\otimes\R^n): A \text{ is symmetric,}\,[A]_{C^{\gamma}}<M_1,\,\norm{A}_{{\infty}}<M_2 \bigr\}.
\end{equation*}
We define
\begin{equation}\label{classA}
\mathcal{A}=\bigl\{A\in C^{\gamma}(\R^n;\R^n\otimes\R^n): \lambda|\xi|^2\leq\langle A(x)\xi,\xi\rangle\leq\Lambda|\xi|^2,\,\forall x,\xi\in\R^n \bigr\}\cap B_{M_1,M_2}.
\end{equation}
\begin{Lem}[Compactness]
\label{Lemma compattezza}
Let $(E_h,u_h)$ be a sequence of $(\kappa_h,\mu)$-minimizers of $\mathcal{F}_{A_h}$ in 
$\Omega$ such that $\sup_h \mathcal F_{A_h}(E_h,u_h;\Omega)<+\infty$, $A_h\rightarrow A_{\infty}$ uniformly on compact sets, where the matrix $A_{\infty}, A_h$ are in the class $\mathcal{A}$ defined in \eqref{classA},
$\kappa_h\rightarrow \kappa\in \mathbb R^+$. There exist a (not relabelled) subsequence and a $(\kappa,\mu)$-minimizer  $(E,u)$ of $\mathcal F_{A_{\infty}}$ in $\Omega$  such that, for every open set $U\Subset \Omega$, it holds
$$
E_h\rightarrow E \mbox { in } L^1(U),\quad u_h\rightarrow u \mbox { in } H^{1}(U),\quad \mathbf{\Phi}_{A_h}(E_h;U)\rightarrow \mathbf{\Phi}_{A_\infty}(E;U).
$$
In addition,
\begin{align}
& \label{boundary1} \mbox{if }x_h\in \partial E_h\cap U \mbox{ and } x_h\rightarrow x \in U, \mbox { then } x\in \partial E \cap U,\\
& \label{boundary2}\mbox{if }x\in \partial E\cap U, \mbox{ there exists } x_h\in \partial E_h\cap U \mbox{ such that } x_h\rightarrow x.
\end{align}
Finally, if we assume also that $\nabla u_h \rightharpoonup 0$ weakly in $L^2_{loc}(\Omega;\mathbb R^n)$ and $\kappa_h \rightarrow 0$, as $h\rightarrow +\infty$, then $E$ is a local minimizer 
of $\mathbf{\Phi}_{A_{\infty}}$, that is
\begin{equation}\label{MinMin}
\mathbf{\Phi}_{A_{\infty}}(E;B_r(x_0))\leq \mathbf{\Phi}_{A_{\infty}}(F;B_r(x_0)),
\end{equation}
for every set $F$ of locally finite perimeter such that $F\vartriangle E \Subset B_r(x_0)\subset \Omega.$
\end{Lem}
\begin{proof}
Using the boundedness condition on $\sup_{h}\mathcal F_{A_h}(E_h,u_h;\Omega)$, we may assume that $u_h$ weakly converges to $u$ in $H^{1}(U)$ and strongly in $L^2(U)$, and 
$\mathbbm{1}_{E_h}$ converges to $\mathbbm{1}_{E}$ in $L^1(U)$, as $h\rightarrow +\infty$. By a lower semicontinuity argument, we start proving the $(\kappa,\mu)$-minimality of $(E,u)$.
Let us fix 
$B_r(x_0)\Subset U$ and assume for simplicity of notation that $x_0=0$. Let $(F,v)$ be a test pair such that $F$ is a set of locally finite perimeter, $F\vartriangle E \Subset B_r$ 
and supp$(u-v)\Subset B_r$. 
Possibly passing to a subsequence and using Fubini's theorem, we may choose $0<r_0<\rho<r$ such that $F\vartriangle E\Subset B_\rho$, $E\setminus B_{r_0}=F\setminus B_{r_0}$, 
 supp$(u-v)\Subset B_{\rho}$, and in addition,
$$
\mathcal{H}^{n-1}(\partial B_{\rho}\cap\partial^*{E})=\mathcal{H}^{n-1}(\partial B_{\rho}\cap \partial^*E_h)=0,
$$
and
\begin{equation}
\label{mismatch}
\lim_{h\rightarrow 0}\mathcal{H}^{n-1}(\partial B_{\rho}\cap(F^{{(1)}}\vartriangle E_h^{{(1)}}))=0.
\end{equation}
Now we choose a cut-off function $\psi\in C_0^1(B_r)$ such that $\psi\equiv 1$ in $B_{\rho}$ and define 
\begin{equation*}
v_h=\psi v+(1-\psi)u_h,\quad\quad 
F_h:=(F\cap B_{\rho})\cup (E_h\setminus B_{\rho})
\end{equation*}
 to test the minimality of $(E_h,u_h)$. Thanks to the $(\kappa_h,\mu)$-minimality of $(E_h,u_h)$ and using also Proposition \ref{compreplace}, we deduce that
\begin{align}\label{min}
& \int_{B_r} \sigma_{E_h}|\nabla u_h|^2 dx+ \mathbf{\Phi}_{A_h}(E_h;B_r)
\leq \int_{B_r} \sigma_{F_h}|\nabla v_h|^2dx+ \mathbf{\Phi}_{A_h}(F_h;B_r)
+\kappa_h 
|F_h\vartriangle E_h|^{\frac{\mu+n-1}{n}}
\notag\\
& \leq \int_{B_r} \sigma_{E_h}(1-\psi)|\nabla u_h|^2 dx+\int_{B_r} \sigma_{F_h}\psi|\nabla v|^2 dx+\int_{B_r} \nabla\psi|u-u_h|^2 dx\notag\\
&+\mathbf{\Phi}_{A_h}(F;B_{\rho})+ \mathbf{\Phi}_{A_h}(E_h;B_r\setminus \overline{B}_{\rho})+\mathbf{\Phi}_{A_h}(B_{\rho};F\vartriangle E_h)+
\kappa_h
|F_h\vartriangle E_h|^{\frac{\mu+n-1}{n}}.
\end{align}
Using the uniform convergence $A_h\rightarrow A_{\infty}$, the strong convergence $u_h\rightarrow u$ in $L^2$, condition \eqref{mismatch}, and getting rid of common terms, from the latter estimate we can write:
\begin{align}
& \int_{B_r} \sigma_{E_h}\psi|\nabla u_h|^2 dx+ \mathbf{\Phi}_{A_{\infty}}(E_h;B_{\rho})
 \notag\\
& \leq \int_{B_r} \sigma_{F_h}\psi|\nabla v|^2 dx
+\mathbf{\Phi}_{A_{\infty}}(F;B_{\rho})
+\kappa_h 
|F_h\vartriangle E_h|^{\frac{\mu+n-1}{n}}
+\varepsilon_h,
\end{align}
for some $\varepsilon_h\rightarrow 0$.
By the lower semicontinuity of the anisotropic perimeter (see \cite[Proposition 3.1]{Simm}),  the equi-integrability of $\left(\nabla u_h\right)_{h\in\N}$ and the lower semicontinuity of Dirichlet integral, we infer that
\begin{equation}
\int_{B_r} \sigma_{E}\psi|\nabla u|^2 dx+ \mathbf{\Phi}_{A_{\infty}}(E;B_{\rho})
\leq \int_{B_r} \sigma_{F}\psi|\nabla v|^2 dx
+\mathbf{\Phi}_{A_{\infty}}(F;B_{\rho})
+\kappa 
|F\vartriangle E|^{\frac{\mu+n-1}{n}}.
\end{equation}
Letting $\psi\downarrow \chi_{B_{\rho}}$ we get
\begin{equation}\label{min1}
\int_{B_{\rho}} \sigma_{E}|\nabla u|^2 dx+ \mathbf{\Phi}_{A_{\infty}}(E;B_{\rho})
\leq \int_{B_{\rho}} \sigma_{F}|\nabla v|^2 dx
+\mathbf{\Phi}_{A_{\infty}}(F;B_{\rho})
+\kappa 
|F\vartriangle E|^{\frac{\mu+n-1}{n}}.
\end{equation}
Similarly, choosing $E=F$ and $u=v$ in \eqref{min}, and arguing as before we get
\begin{equation}
\limsup_{h\rightarrow +\infty}\bigg(\int_{B_{\rho}} \sigma_{E_h}\psi|\nabla u_h|^2 dx+ \mathbf{\Phi}_{A_{\infty}}(E_h;B_{\rho})\bigg)
\leq \int_{B_{\rho}} \sigma_{E}\psi|\nabla u|^2 dx+ \mathbf{\Phi}_{A_{h}}(E;B_{\rho}).
\end{equation}
Letting $\psi\downarrow \chi_{B_{\rho}}$ we conclude
\begin{equation}
\lim_{h\rightarrow +\infty}\mathbf{\Phi}_{A_{h}}(E_h;B_{\rho})=\mathbf{\Phi}_{A_{\infty}}(E;B_{\rho}),\quad \quad \quad 
\lim_{h\rightarrow +\infty}\int_{B_{\rho}} \sigma_{E_h}|\nabla u_h|^2 dx=\int_{B_{\rho}} \sigma_{E}|\nabla u|^2 dx.
\end{equation}
With a usual argument we can deduce
$u_h\rightarrow u$ in $W^{1,2}(U)$ and $\mathbf{\Phi}_{A_{h}}(E_h;U)\rightarrow \mathbf{\Phi}_{A_{\infty}}(E;U)$, for every open set $U\Subset \Omega$.
The topological information stated in \eqref{boundary1} and \eqref{boundary2} follows as in \cite[Theorem 21.14]{Ma},
indeed they are a consequence of the lower and upper density estimates given above. Finally, if $\nabla u_h \rightharpoonup 0$ weakly in $L^2_{loc}(\Omega;\mathbb R^n)$ and $\kappa_h \rightarrow 0$, we can choose $v=u$ in \eqref{min1}, deriving \eqref{MinMin}.
\end{proof}

\section{Reverse Poincaré inequality}
\label{Reverse Poincaré inequality}
In this section we derive a reverse Poincaré inequality which lets us estimate the excess around a point of the boundary of the transformed set with its flatness. The first step in the proof is to establish a weak form of this inequality.\\
\indent In the following proposition, it is proved that if the anisotropy matrix valued in a point $x_0$ is the identity, then around $x_0$ the anisotropic perimeter is comparable to the perimeter.

\begin{Prop}\label{MimimalitàCasoParticolare}
Let $x_0\in\Omega$ and $r>0$. There exists a positive constant $c_9=\break c_9\big(n,\alpha,\beta,\lambda,\Lambda,\kappa,\mu,\norm{\D u}_{L^2(\Omega)}\big)$ such that if $(E,u)$ is a $(\kappa r^\mu,\mu)$-minimizer of $\mathcal{F}_{x_0,A_{x_0}\circ \Psi_{x_0,r}^{-1}}$ in $B_1$, with $0\in \dd E\cap B_1$, then 
    \begin{align*}
&\int_{B_\rho}\sigma_{E}|\D uA^{-\frac{1}{2}}(x_0)|^2\,dx+P(E;B_\rho)\\
&\leq\int_{B_\rho}\sigma_{ F}|\D  vA^{-\frac{1}{2}}(x_0)|^2\,dx+P( F;B_\rho)+c_9\big(\kappa+[A]_{C^\gamma}\big)r^\mu\rho^{n-1+\mu},
    \end{align*}
    for every $(F,v)$ such that $F\vartriangle E\Subset B_\rho\subset B_1$ and $v\in u+H^1_0(B_\rho)$.
\end{Prop}
\begin{proof}
    Let $(F,v)$ be such that $F\vartriangle E\Subset B_\rho$ and $v\in u+H^1_0(B_\rho)$. We can assume that
    \begin{align*}
        \int_{B_\rho}\sigma_{E}|\D u A^{-\frac{1}{2}}(x_0)|^2\,dx+P(E;B_\rho)\geq\int_{B_\rho}\sigma_{F}|\D v A^{-\frac{1}{2}}(x_0)|^2\,dx+P(F;B_\rho).
    \end{align*}
    We remark that $A_{x_0}\circ \Psi_{x_0,r}^{-1}$ is H\"older continuous and 
    \begin{equation*}
     \big[A_{x_0}\circ \Psi_{x_0,r}^{-1}\big]_{C^\mu}\leq \frac{\Lambda^\frac{\mu}{2}}{\lambda}[A]_{C^\mu}
     r^{\mu}.
    \end{equation*}
    Since $\big(A_{x_0}\circ \Psi_{x_0,r}^{-1}\big)(0)=I$, by the H\"older continuity of $A_{x_0}\circ \Psi_{x_0,r}^{-1}$ we infer\\
    \begin{align*}
    |\nu_{E}|
    & =\prodscal{\big(A_{x_0}\circ \Psi_{x_0,r}^{-1}\big)(0)\;\nu_{ E}}{\nu_{E}}^{\frac{1}{2}}\leq \prodscal{\big(A_{x_0}\circ \Psi_{x_0,r}^{-1}\big)(x)\;\nu_{E}}{\nu_{ E}}^{\frac{1}{2}}+\frac{1}{2\lambda}\big[A_{x_0}\circ \Psi_{x_0,r}^{-1}\big]_{C^\mu}\rho^\mu\\
    & \leq \prodscal{\big(A_{x_0}\circ \Psi_{x_0,r}^{-1}\big)(x)\;\nu_{ E}}{\nu_{ E}}^{\frac{1}{2}}+\frac{\Lambda^{\frac{\mu}{2}}}{2\lambda^2}[A]_{C^\mu}(r\rho)^\mu,
    \end{align*}
    for any $x\in B_\rho$. Integrating over $B_\rho$ with respect to the measure $\mathcal{H}^{n-1}\resmes\dd^*{E}$ and adding to both sides the term $\int_{B_\rho}\sigma_{ E}|\D u A^{-\frac{1}{2}}(x_0)|^2$, we obtain
\begin{align*}
\int_{B_\rho}\sigma_{ E}|\D u A^{-\frac{1}{2}}(x_0)|^2\,dx+P(E;B_\rho)
& \leq\int_{B_\rho} \sigma_{E}|\D u A^{-\frac{1}{2}}(x_0)|^2\,dx+\mathbf{\Phi}_{A_{x_0}\circ \Psi_{x_0,r}^{-1}}(E;B_\rho)\\
& +\frac{\Lambda^{\frac{\mu}{2}}}{2\lambda^2}[A]_{C^\mu}(r\rho)^\mu P({E};B_\rho).
    \end{align*}
    Arguing in a similar way, we get
    \begin{align*}
\int_{B_\rho}\sigma_{F}|\D v A^{-\frac{1}{2}}(x_0)|^2\,dx+\mathbf{\Phi}_{A_{x_0}\circ \Psi_{x_0,r}^{-1}}(F;B_\rho)
& \leq\int_{B_\rho}\sigma_{F}|\D v A^{-\frac{1}{2}}(x_0)|^2\,dx+P(F;B_\rho)\\
& +\frac{\Lambda^{\frac{\mu}{2}}}{2\lambda^2}[A]_{C^\mu}(r\rho)^\mu P(F;B_\rho).
    \end{align*}
Applying the definition of $(\kappa r^\mu,\mu)$-minimality of $(E,u)$ and using the previous two  inequalities, we write
    \begin{align*}
       & \int_{B_\rho}\sigma_{E}|\D u A^{-\frac{1}{2}}(x_0)|^2\,dx+P(E;B_\rho)\\
       & \leq\int_{B_\rho}\sigma_{E}|\D u A^{-\frac{1}{2}}(x_0)|^2\,dx+\mathbf{\Phi}_{A_{x_0}\circ \Psi_{x_0,r}^{-1}}(E;B_\rho)+\frac{\Lambda^{\frac{\mu}{2}}}{2\lambda^2}[A]_{C^\mu}(r\rho)^\mu P(E;B_\rho)\\
       & \leq\int_{B_\rho}\sigma_{F}|\D v A^{-\frac{1}{2}}(x_0)|^2\,dx+\mathbf{\Phi}_{A_{x_0}\circ \Psi_{x_0,r}^{-1}}(F;B_\rho)+\kappa r^\mu |E\vartriangle F|^{\frac{n-1+\mu}{n}}+\frac{\Lambda^{\frac{\mu}{2}}}{2\lambda^2}[A]_{C^\mu}(r\rho)^\mu P(E;B_\rho)\\
       & \leq\int_{B_\rho}\sigma_{F}|\D v A^{-\frac{1}{2}}(x_0)|^2\,dx+P(F;B_\rho)+c(n)\kappa r^\mu \rho^{n-1+\mu}\\
       & +\frac{\Lambda^{\frac{\mu}{2}}}{2\lambda^2}[A]_{C^\mu}(r\rho)^\mu [P(E;B_\rho)+P(F;B_\rho)]\\
       & \leq\int_{B_\rho}\sigma_{F}|\D v A^{-\frac{1}{2}}(x_0)|^2\,dx+P(F;B_\rho)+c(n)\kappa r^\mu\rho^{n-1+\mu}\\
       & +\frac{\Lambda^{\frac{\mu}{2}}}{2\lambda^2}[A]_{C^\mu}(r\rho)^\mu \bigg[2P(E;B_\rho)+\int_{B_\rho}\sigma_{E}|\D u|^2\,dx\bigg]\\
       & \leq\int_{B_\rho}\sigma_{F}|\D v A^{-\frac{1}{2}}(x_0)|^2\,dx+P(F;B_\rho)+c(n,\lambda,\Lambda,c_3)\Big(\kappa+[A]_{C^\gamma}\Big)r^\mu\rho^{n-1+\mu},
    \end{align*}
    where $c_3$ is the constant appearing in Theorem \ref{Energy upper bound}, which leads to the thesis.
\end{proof}

At this point, we are able to establish a weak form of the reverse Poincaré inequality. The strategy for its proof is the same outlined in \cite[Lemma 24.9]{Ma} (see also \cite[Lemma 7.3]{Simm} or \cite[Lemma 10]{EL2}).
\begin{Lem}[Weak reverse Poincaré inequality]
\label{Weak reverse Poincaré inequality}
\label{Weak Reverse Poincaré}
Let $x_0\in\Omega$ and $r>0$. If $(E,u)$ is a $(\kappa r^{\mu},\mu)$-minimizer of $\mathcal{F}_{x_0,A_{x_0}\circ\Psi_{x_0,r}^{-1}}$ in $\mathbf{C}_4$ such that
\begin{equation*}
|x_n|<\frac{1}{8},\quad\forall x\in \mathbf{C}_2\cap\dd E,
\end{equation*}
\begin{equation*}
\left|\left\{ x\in \mathbf{C}_2\setminus E\,:\, x_n<-\frac{1}{8} \right\}\right|=\left|\left\{ x\in \mathbf{C}_2\cap E\,:\, x_n>\frac{1}{8} \right\}\right|=0,
\end{equation*}
and if $z\in\R^{n-1}$ and $s>0$ are such that
\begin{equation}
\label{eqqqz18}
\mathbf{K}_s(z)\subset \mathbf{C}_2,\qquad \mathcal{H}^{n-1}(\dd E\cap\dd \mathbf{K}_s(z))=0,
\end{equation}
then, for every $|c|<\frac{1}{4}$,
\begin{align*}
& P(E;\mathbf{K}_{\frac{s}{2}}(z))-\mathcal{H}^{n-1}(\mathbf{D}_{\frac{s}{2}}(z))\leq c_{10}\Bigg\{\bigg[\left(P(E;\mathbf{K}_s(z))-\mathcal{H}^{n-1}(\mathbf{D}_s(z))\right)\\
& \times\int_{\mathbf{K}_s(z)\cap\dd^*E}\frac{(x_n-c)^2}{s^2}d\,\mathcal{H}^{n-1}\bigg]^{\frac{1}{2}}+\int_{\mathbf{K}_{s}(z)}|\D u|^2\,dx+\bigl(\kappa+[A]_{C^\mu}\bigr)r^\mu\Bigg\},
\end{align*}
for some positive constant $c_{10}=c_{10}\big(n,\alpha,\beta,\lambda,\Lambda,\kappa,\mu,\norm{\D u}_{L^2(\Omega)}\big)$.
\end{Lem}
\begin{proof}
We may assume that $z=0$. The set function
\begin{equation*}
m(G)=P(E;\mathbf{C}_2\cap p^{-1}(G))-\mathcal{H}^{n-1}(G),\quad\text{for }G\subset \mathbf{D}_2,
\end{equation*}
defines a Radon measure on $\R^{n-1}$, supported in $\mathbf{D}_2$. Since $E$ is a set of locally finite perimeter, by \cite[Theorem 13.8]{Ma} there exist a sequence $(E_h)_{h\in\N}$ of open subsets of $\R^n$ with smooth boundary and a vanishing sequence $(\varepsilon_h)_{h\in\N}\subset\R^+$ such that
\begin{equation*}
E_h\overset{loc}{\rightarrow}E,\quad\mathcal{H}^{n-1}\resmes\dd E_h\rightarrow\mathcal{H}^{n-1}\resmes\dd E,\quad \dd E_h\subset I_{\varepsilon_h}(\dd E),
\end{equation*}
as $h\rightarrow+\infty$, where $I_{\varepsilon_h}(\dd E)$ is a tubular neighborhood of $\dd E$ with half-lenght $\varepsilon_h$. By the coarea formula we get
\begin{equation*}
\mathcal{H}^{n-1}(\dd \mathbf{K}_{\rho s}\cap(E^{(1)}\vartriangle E_h))\rightarrow0,\quad \text{for a.e. }\rho\in\left(\frac{2}{3},\frac{3}{4}\right).
\end{equation*}
Moreover, provided $h$ is large enough, by $\dd E_h\subset I_{\varepsilon_h}(\dd E)$, we get:
\begin{equation*}
|x_n|<\frac{1}{4},\quad\forall x\in \mathbf{C}_2\cap\dd E_h,
\end{equation*}
\begin{equation*}
\left\{ x\in \mathbf{C}_2\,:\, x_n<-\frac{1}{4} \right\}\subset \mathbf{C}_2\cap E_h\subset\left\{ x\in \mathbf{C}_2\,:\, x_n<\frac{1}{4} \right\}.
\end{equation*}
Therefore, given $\lambda\in\left(0,\frac{1}{4}\right)$ and $|c|<\frac{1}{4}$, we are in position to apply \cite[Lemma 24.8]{Ma} to every $E_h$ to deduce that there exists $I_h\subset\left( \frac{2}{3},\frac{3}{4} \right)$, with $|I_h|\geq\frac{1}{24}$, and, for any $\rho\in I_h$, there exists an open subset $F_h$ of $\R^n$ of locally finite perimeter such that
\begin{equation}
\label{eqqq15}
F_h\cap\dd \mathbf{K}_{\rho s}=E_h\cap\dd \mathbf{K}_{\rho s},\\
\end{equation} 
\begin{equation}
\label{eqqq16}
\mathbf{K}_{\frac{s}{2}}\cap\dd F_h=\mathbf{D}_{\frac{s}{2}}\times\{c\},
\end{equation}
\begin{align}
\label{eqqq17}
P(F_h;\mathbf{K}_{\rho s})-\mathcal{H}^{n-1}(\mathbf{D}_{\rho s})\leq 
& c(n)\bigg\{ \lambda\left( P(E_h;\mathbf{K}_s)-\mathcal{H}^{n-1}(\mathbf{D}_s) \right)+\frac{1}{\lambda}\int_{\mathbf{K}_s\cap\dd E_h}\frac{|x_n-c|^2}{s^2}\,d\mathcal{H}^{n-1} \bigg\}.
\end{align}
Clearly $\displaystyle\bigcap_{h\in\N}\bigcup_{k\geq h}|I_k|\geq\frac{1}{24}>0$ and thus there exist a divergent subsequence $(h_k)_{k\in\N}$ and $\rho\in\left( \frac{2}{3},\frac{3}{4}\right)$ such that
\begin{equation*}
\rho\in\bigcap_{k\in\N} I_{h_k} \quad\text{and}\quad \lim_{k\rightarrow+\infty}\mathcal{H}^{n-1}(\dd \mathbf{K}_{\rho s}\cap(E^{(1)}\vartriangle E_{h_k}))=0.
\end{equation*}
We will write $F_k$ in place of $F_{h_k}$. We consider the comparison set $G_k=(F_k\cap \mathbf{K}_{\rho s})\cup (E\setminus \mathbf{K}_{\rho  s})$. By applying \cite[formula (16.33)]{Ma} we infer that
\begin{equation*}
P(G_k;\mathbf{K}_s)=P(F_k;\mathbf{K}_{\rho s})+P(E;\mathbf{K}_s\setminus \mathbf{K}_{\rho s})+\sigma_k,
\end{equation*}
where, thanks to \eqref{eqqq15}, $\sigma_k=\mathcal{H}^{n-1}(\dd \mathbf{K}_{\rho s}\cap (E^{(1)}\vartriangle F_k))=\mathcal{H}^{n-1}(\dd \mathbf{K}_{\rho s}\cap (E^{(1)}\vartriangle E_{h_k}))\rightarrow0$,  as $k\rightarrow+\infty$. We apply Proposition \ref{MimimalitàCasoParticolare}, deducing the following relation:
    \begin{align*}
&\int_{\mathrm{spt}(u-v)}\sigma_{E}|\D u A^{-\frac{1}{2}}(x_0)|^2\,dx+P(E;B_{\tilde{\rho}})\\
&\leq\int_{\mathrm{spt}(u-v)}\sigma_G|\D vA^{-\frac{1}{2}}(x_0)|^2\,dx+P(G;B_{\tilde{\rho}})+c\big(n,\alpha,\beta,\lambda,\Lambda,\kappa,\mu,\norm{\D u}_{L^2(\Omega)}\big)\big(\kappa r^{\mu}+[A]_{C^\mu}r^\mu\big)\tilde{\rho}^{n-1+\mu},
    \end{align*}
    for every $(G,v)$ such that $G\vartriangle E\Subset B_{\tilde\rho}\subset \mathbf{C}_4$ and $v\in u+H^1_0(B_{\tilde{\rho}})$.

Now we test the previous relation of minimality with $(G_k,u)$, as $E\vartriangle G_k\Subset \mathbf{K}_s\subset B_4\subset \mathbf{C}_4$, and get rid of the common terms obtaining
\begin{align}
P(E;\mathbf{K}_{\rho s})
& \leq P(F_k;\mathbf{K}_{\rho s})+\sigma_k+c\big(n,\alpha,\beta,\lambda,\Lambda,\kappa,\mu,\norm{\D u}_{L^2(\Omega)}\big)\bigg[\int_{\mathbf{K}_{\rho s}}|\D u|^2\,dx+\big(\kappa r^{\mu}+[A]_{C^\mu} r^\mu\big)\bigg].
\label{eqqqqq1}
\end{align}
 Thus, since $m$ is nondecreasing and $\rho\in\big(\frac{2}{3},\frac{3}{4}\big)$, by \eqref{eqqqqq1} and \eqref{eqqq17} we deduce that
\begin{align*}
& P(E;\mathbf{K}_{\frac{s}{2}})-\mathcal{H}^{n-1}(\mathbf{D}_{\frac{s}{2}})=m(\mathbf{D}_{\frac{s}{2}})\leq m(\mathbf{D}_{\rho s})=P(E;\mathbf{K}_{\rho s})-\mathcal{H}^{n-1}(\mathbf{D}_{\rho s})\\
& \leq P(F_k;\mathbf{K}_{\rho s})-\mathcal{H}^{n-1}(\mathbf{D}_{\rho s})+\sigma_k+c\bigg[\int_{\mathbf{K}_{\rho s}}|\D u|^2\,dx+\big(\kappa +[A]_{C^\mu}\big)r^\mu\bigg]\\
& \leq c(n)\bigg\{ \lambda\left( P(E_{h_k};\mathbf{K}_s)-\mathcal{H}^{n-1}(\mathbf{D}_s) \right)+\frac{1}{\lambda}\int_{\mathbf{K}_s\cap\dd E_{h_k}}\frac{|x_n-c|^2}{s^2}\,d\mathcal{H}^{n-1} \bigg\}\\
& +c \bigg[\int_{\mathbf{K}_{s}}|\D u|^2\,dx+\bigl(\kappa+[A]_{C^\mu}\bigr)r^\mu\bigg],
\end{align*}
where $c=c\big(n,\alpha,\beta,\lambda,\Lambda,\kappa,\mu,\norm{\D u}_{L^2(\Omega)}\big)$.
Letting $k\rightarrow +\infty$, \eqref{eqqqz18} implies that $P(E_{h(k)};\mathbf{K}_s)\rightarrow P(E;\mathbf{K}_s)$ and therefore
\begin{align}
\label{eqqq19}
P(E;\mathbf{K}_{\frac{s}{2}})-\mathcal{H}^{n-1}(\mathbf{D}_{\frac{s}{2}})
& \leq c \bigg\{ \lambda\left( P(E;\mathbf{K}_s)-\mathcal{H}^{n-1}(\mathbf{D}_s) \right)+\frac{1}{\lambda}\int_{\mathbf{K}_s\cap\dd E}\frac{|x_n-c|^2}{s^2}\,d\mathcal{H}^{n-1}\notag\\
& +\int_{\mathbf{K}_{rs}}|\D u|^2\,dx+\bigl(\kappa+[A]_{C^\mu}\bigr)r^\mu\bigg\},
\end{align}
for any $\lambda\in\left(0,\frac{1}{4}\right)$. If $\lambda>\frac{1}{4}$, then
\begin{align*}
& P(E;\mathbf{K}_{\frac{s}{2}})-\mathcal{H}^{n-1}(\mathbf{D}_{\frac{s}{2}})=m(\mathbf{D}_{\frac{s}{2}})\leq m(\mathbf{D}_{\rho s})\\
& \leq 4\lambda P(E;\mathbf{K}_{\rho s})-\mathcal{H}^{n-1}(\mathbf{D}_{\rho s})\leq c(n)\lambda\left( P(E;\mathbf{K}_s)-\mathcal{H}^{n-1}(\mathbf{D}_s) \right),
\end{align*}
and thus \eqref{eqqq19} holds  true for $\lambda>0$, provided we choose $c(n)\geq 4$. Minimizing over $\lambda$, we get the thesis.
\end{proof}

Finally, we are able to prove the main result of this section.

\begin{Thm}[Reverse Poincaré inequality]
\label{RevPI}
Let $x_0\in\Omega$ and $r>0$. There exist two positive constants $c_{11}=c_{11}\big(n,\alpha,\beta,\lambda,\Lambda,\kappa,\mu,\norm{\D u}_{L^2(\Omega)}\big)$ and $\varepsilon_5=\varepsilon_5(n)$ such that if $(E,u)$ is a $(\kappa r^{\mu},\mu)$-minimizer of $\mathcal{F}_{x_0,A_{x_0}\circ\Psi_{x_0,r}^{-1}}$ in $\mathbf{C}_{4\tau}(0,\nu)$, with $0\in \partial E$, $\tau >0$ and
\begin{equation*}
\ecc(0,4\tau,\nu)<\varepsilon_5,
\end{equation*}
then
\begin{align}
\label{RevPoinc}
\ecc(0,\tau,\nu)
\leq
c_{11}\biggl(\frac{1}{\tau^{n+1}}\int_{\partial E\cap \mathbf{C}_{2 \tau}(0,\nu)}
&|\left<\nu,x\right>-c|^2d\mathcal{H}^{n-1}\\
&+\frac{1}{\tau^{n-1}}\int_{\mathbf{C}_{2\tau}(0,\nu)}|\D u|^2\,dx+ \bigl(\kappa+[A]_{C^\mu}\bigr)(\tau r)^\mu\biggr),\notag
\end{align}
for every $c\in \mathbbm R$.
\end{Thm}

\begin{proof}
The proof of this result follows the same strategy employed in \cite[Theorem 6]{EL2}. We emphasize only small differences between the two proofs. Up to a rotation and employing a usual scaling argument, by Proposition \ref{PropScaling}, with a small abuse of notation, we may assume that $(E,u)$ is a $(\kappa (\tau r)^\mu,\mu)$-minimizer of $\mathcal{F}_{x_0,A_{x_0}\circ\Psi_{x_0,\tau r}^{-1}}$ in $\mathbf{C}_4$, with $0\in\dd \tilde{E}$. Leveraging the compactness of the perimeter and Theorem \ref{LipApp}, it is possible to show that 
\begin{equation*}
|x_n|<\frac{1}{4},\quad\forall x\in \mathbf{C}_2\cap\dd E,
\end{equation*}
\begin{equation*}
\left|\left\{ x\in \mathbf{C}_2\setminus E\,:\, x_n<-\frac{1}{8} \right\}\right|=\left|\left\{ x\in \mathbf{C}_2\cap E\,:\, x_n>\frac{1}{8} \right\}\right|=0.
\end{equation*}
Thus, for any $z\in\R^{n-1}$ and $s>0$ such that
\begin{equation}
\label{eqqq18}
\mathbf{K}_s(z)\subset \mathbf{C}_2,\quad \mathcal{H}^{n-1}(\dd E\cap\dd \mathbf{K}_s(z))=0,
\end{equation}
we apply Lemma \ref{Weak reverse Poincaré inequality}, deducing that, for every $|c|<\frac{1}{4}$,
\begin{align}
\label{WRPI} 
P(E;\mathbf{K}_{s}(z))-\mathcal{H}^{n-1}(\mathbf{D}_s(z))
& \leq c\left\{ \bigg[[P(E;\mathbf{K}_{2s}(z))-\mathcal{H}^{n-1}(\mathbf{D}_{2s}(z))]\inf_{|c|<\frac{1}{4}}\int_{\mathbf{C}_2\cap\dd E}|x_n-c|^2\,d\mathcal{H}^{n-1}\bigg]^{\frac{1}{2}}\right.\notag\\
& \left.+ \int_{\mathbf{K}_{s}}|\D u|^2\,dx+\kappa
\tau r^\mu+[A]_{C^\mu}(\tau r)^\mu \right\},
\end{align}
for some positive constant $c=c\big(n,\alpha,\beta,\lambda,\Lambda,\kappa,\mu,\norm{\D u}_{L^2(\Omega)}\big)$. Hence, proceeding as in \cite[Theorem 6]{EL2}, by a covering argument, it is possible to show that \eqref{WRPI} implies \eqref{RevPoinc}.
\end{proof}

\section{Proof of the main theorem}
\label{Proof of the main theorem}

The strategy adopted to establish the main result involves two key steps: first proving a first variation formula for the bulk energy of $\mathcal{F}_{x_0,A_{x_0}\circ\Psi_{x_0,r}^{-1}}$, then establishing an excess improvement theorem for transformed couples, which in turn implies an analogous theorem for the original ones.

\begin{Prop}[First variation formula for the bulk term]\label{VarB}
$x_0\in\Omega$, $u\in H^1(B_1)$ and $X\in C^1_c(B_1;\R^n)$. We define $\Phi_t(x)=x+tX(x)$, for any $x\in\R^n$ and $t>0$. Accordingly, we define
\begin{equation}
E_t:=\Phi_t(E),\quad u_t:=u\circ\Phi_t^{-1}.
\end{equation}
There exist two constants $c_{12}=c_{12}(\beta,\lambda,\D X)>0$ and $t_0>0$ such that it holds that
\begin{equation*}
\int_{B_1}\sigma_{E_t}|\D u_tA^{-\frac{1}{2}}(x_0)|^2\,dx-\int_{B_1}\sigma_E|\D uA^{-\frac{1}{2}}(x_0)|^2\,dx\leq c_{12}(t+o(t))\int_{B_1}|\D u|^2\,dx,
\end{equation*}
for any $0<t<t_0$.
\end{Prop}
\begin{proof}
Taking into account that
\begin{equation*}
\D \Phi_t^{-1}(\Phi_t(x))=I-t\D X(x)+o(t), \quad \textnormal{J}\Phi_t(x)=1+t\textnormal{div}X(x)+o(t),
\end{equation*}
for any $x\in\R^n$ and $t>0$, by the change of variable $y=\Phi_t(x)$ we obtain
\begin{align}\label{EL1}
& \int_{B_1}\sigma_{E_t}|\D u_tA^{-\frac{1}{2}}(x_0)|^2\,dy-\int_{B_1}\sigma_E|\D uA^{-\frac{1}{2}}(x_0)|^2\,dx\notag\\
= & \int_{B_1}\sigma_E\big|[\D u-t\D u\D X+\D u \,o(t)]A^{-\frac{1}{2}}(x_0)\big|^2(1+t\mathrm{div}X+o(t))\,dx- \int_{B_1}\sigma_E|\D uA^{-\frac{1}{2}}(x_0)|^2\,dx\notag\\
& =\int_{B_1}\sigma_E\big[|\D uA^{-\frac{1}{2}}(x_0)|^2+|\D uA^{-\frac{1}{2}}(x_0)|^2(t\mathrm{div}X+o(t))\big]\,dx+H(t,\D u,\D X)\notag\\
& - \int_{B_1}\sigma_E|\D uA^{-\frac{1}{2}}(x_0)|^2\,dx\notag\\
& = \int_{B_1}\sigma_E|\D uA^{-\frac{1}{2}}(x_0)|^2(t\mathrm{div}X+o(t))\,dx+H(t,\D u,\D X),
\end{align}
where
\begin{align*}
& H(t,\D u,\D X)
= \int_{B_1}\sigma_E\big|[-t\D u\D X+\D u \,o(t)]A^{-\frac{1}{2}}(x_0)\big|^2(1+t\mathrm{div}X+o(t))\,dx\\
& + \int_{B_1}2\prodscal{\D u A^{-\frac{1}{2}}(x_0)}{(-t\D u\D X+\D u\,o(t))A^{-\frac{1}{2}}(x_0)}(1+t\mathrm{div}X+o(t))\,dx.
\end{align*}
We estimate
\begin{equation}\label{EL2}
\int_{B_1}\sigma_E|\D uA^{-\frac{1}{2}}(x_0)|^2(t\mathrm{div}X+o(t))\,dx\leq c(\beta,\lambda,\D X)(t+o(t))\int_{B_1}|\D u|^2\,dx
\end{equation}
and
\begin{align}\label{EL3}
H(t,\D u,\D X)
& \leq c(\beta,\lambda,\D X)\int_{B_1}(t+o(t))^2|\D u|^2(1+t\mathrm{div}X+o(t))\,dx\notag\\
& +c(\beta,\lambda,\D X)\int_{B_1}(t+o(t))|\D u|^2(1+t\mathrm{div}X+o(t))\,dx\notag\\
& \leq c(\beta,\lambda,\D X)(t+o(t))\int_{B_1}|\D u|^2\,dx.
\end{align}
Inserting $\eqref{EL3}$ and $\eqref{EL2}$ in $\eqref{EL1}$ we get the desired inequality.
\end{proof}

Here we present the proof of the excess improvement theorem for transformed couples.
\begin{Thm}[Excess improvement for the transformed couple]\label{Teorema eccesso raddrizzato}
For any $\omega\in(0,1)$, $ \tilde{\sigma}\in(0,1)$, $\tilde{M}>0$, $\tilde{\tau}\in\big(0,\frac{1}{16}\big)$ there exists a constant $\tilde{\varepsilon}=\tilde{\varepsilon}(\tilde{\sigma},\tilde{M},\tilde{\tau})>0$ such that if $(\tilde{E},\tilde{u})$ is a $(\tilde{\kappa},\mu)$-minimizer of $\mathcal{F}_{x_0,A_{x_0}}$ in $B_{\tilde r}(x_0)$, with $x_0\in\dd\tilde{E}$, such that
\begin{equation*}
\mathbf{e}(\tilde{E},x_0,\tilde{r})\leq\tilde{\varepsilon},\quad\mathcal{D}_{\tilde{u}}(x_0,\tilde{r})+\tilde{r}^{(1-\omega)\mu}\leq \tilde{M}\mathbf{e}(\tilde{E},x_0,\tilde{\sigma}\tilde{r}),
\end{equation*}
then there exists a constant $c_{13}=c_{13}\big(n,\alpha,\beta,\lambda,\Lambda,\kappa,\mu,\norm{\D u}_{L^2(\Omega)}\big)>0$ such that
\begin{equation*}
\mathbf{e}(\tilde{E},x_0,\tilde{\tau}\tilde{r})\leq c_{13}\big(\tilde{\tau}^2\mathbf{e}(\tilde{E},x_0,\tilde{r})+\mathcal{D}_{\tilde{u}}(x_0,4\tilde{\tau}\tilde{r})+(\tilde{\tau}\tilde{r})^\mu\big).
\end{equation*}
\end{Thm}
\begin{proof}
Let us assume by contradiction that there exist a vanishing sequence $(\tilde{r}_h)_{h\in\N}\subset\R^+$ and a sequence $((\tilde{E}_h,\tilde{u}_h))_{h\in\N}$ of $(\tilde{\kappa},\mu)$-minimizers of $\mathcal{F}_{x_0,A}$ in $B_{\tilde{r}_h}(x_0)$, with $x_0\in\dd \tilde{E}_h$, such that
\begin{equation*}
\mathbf{e}(\tilde{E}_h,x_0,\tilde{r}_h)=:\varepsilon_h\rightarrow 0,\quad \mathcal{D}_{\tilde{u}_h}(x_0,\tilde{r}_h)+\tilde{r}_h^{(1-\omega)\mu}\leq \tilde{M}\mathbf{e}(\tilde{E}_h,x_0,\tilde{\sigma}\tilde{r_h}),
\end{equation*}
and
\begin{equation}
\mathbf{e}(\tilde{E}_h,x_0,\tilde{\tau}\tilde{r}_h)> \overline{C}\big(\tilde{\tau}^2\mathbf{e}(\tilde{E}_h,x_0,\tilde{r}_h)+\mathcal{D}_{\tilde{u}_h}(x_0,4\tilde{\tau}\tilde{r}_h)+(\tilde{\tau}\tilde{r}_h)^\mu\big),
\end{equation}
for some constant $\overline{C}>0$ to be chosen. Employing the usual scaling argument and applying Proposition \ref{PropScaling}, with a small abuse of notation we may assume that $((\tilde{E}_h,\tilde{u}_h))_{h\in\N}$ is a sequence of $(\tilde{\kappa}\tilde{r}_h^\mu,\mu)$-minimizers of $\mathcal{F}_{x_0,A_{x_0}\circ\Psi_{x_0,\tilde{r}_h}^{-1}}$ in $B_1$, with $0\in\dd \tilde{E}_h$, such that
\begin{equation}\label{EccDirRag}
\mathbf{e}(\tilde{E}_h,0,1)=\varepsilon_h\rightarrow 0,\quad \mathcal{D}_{\tilde{u}_h}(0,1)+\tilde{r}_h^{(1-\omega)\mu}\leq \tilde{M}\mathbf{e}(\tilde{E}_h,0,\tilde{\sigma}),
\end{equation}
and
\begin{equation}
\mathbf{e}(\tilde{E}_h,0,\tilde{\tau})> \overline{C}\big(\tilde{\tau}^2\mathbf{e}(\tilde{E}_h,0,1)+\mathcal{D}_{\tilde{u}_h}(0,4\tilde{\tau})+(\tilde{\tau}\tilde{r}_h)^\mu\big).
\end{equation}
Up to rotating each $\tilde{E}_h$ we may also assume that, for all $h\in\N$,
\begin{equation*}
{\mathbf e}(\tilde{E}_h,0,1)=\frac{1}{2}\int_{\dd \tilde{E}_h\cap B_1}\abs{\nu_{E_h}-e_n}^2\, d\mathcal{H}^{n-1}.
\end{equation*}
\textbf{Step 1.} Thanks to the Lipschitz approximation theorem, for $h$ sufficiently large, there exists a 1-Lipschitz function $f_h\colon\R^{n-1}\rightarrow\R$ such that
\begin{equation}
\label{1}
\sup_{\R^{n-1}}\abs{f_h}\leq c_8\varepsilon_h^{\frac{1}{2(n-1)}}, \quad \mathcal{H}^{n-1}((\dd \tilde{E}_h\vartriangle\Gamma_{f_h})\cap B_{\frac{1}{2}})\leq c_8\varepsilon_h, \quad \int_{ {\mathbf{D}_{\frac{1}{2}}}}\abs{\D' f_h}^2\,dx'\leq c_8\varepsilon_h.
\end{equation}
We define
\begin{equation*}
g_h(x'):=\frac{f_h(x')-a_h}{\sqrt{\varepsilon_h}}, \quad\text{where}\quad a_h=\dashint_{ {\mathbf{D}_{\frac{1}{2}}}}f_h\,dx',
\end{equation*}
and we assume, up to a subsequence, that $\lbrace g_h \rbrace_{h\in\N}$ converges weakly in $H^1( {\mathbf{D}_{\frac{1}{2}}})$ and strongly in $L^2( {\mathbf{D}_{\frac{1}{2}}})$ to a function $g$. We prove that $g$ is harmonic in $ {\mathbf{D}_{\frac{1}{2}}}$. It is enough to show that
\begin{equation}
\label{2}
\lim_{h\rightarrow +\infty}\frac{1}{\sqrt{\varepsilon_h}}\int_{ {\mathbf{D}_{\frac{1}{2}}}}\frac{\langle \D' f_h,\D'\phi\rangle}{\sqrt{1+\abs{\D' f_h}^2}}\,dx'=0,
\end{equation}
for all $\phi\in C_0^1( {\mathbf{D}_{\frac{1}{2}}})$. We fix $\delta>0$ so that supp$\,\phi\times [-2\delta,2\delta]\subset B_{\frac{1}{2}}$ and choose a cut-off function $\psi\colon\R\rightarrow[0,1]$ with supp$\,\psi\subset (-2\delta,2\delta)$, $\psi=1$ in $(-\delta,\delta)$. Let us define
$$
\Phi_{h}(x):=x+\tilde{r}_h^{\omega\mu} X(x),\quad \text{where }X(x)=\phi(x')\psi(x_n) e_n,
$$
for $x\in\R^n$. We apply Proposition \ref{MimimalitàCasoParticolare} to deduce that
\begin{align}\label{FV0}
P(\tilde{E}_h;B_{\frac{1}{2}})-P(\Phi_{h}(\tilde{E}_h);B_{\frac{1}{2}})
& \leq \int_{B_{\frac{1}{2}}}\sigma_{\Phi_{h}(\tilde{E}_h)}|\D (\tilde{u}_h\circ\Phi^{-1}_{h})A^{-\frac{1}{2}}(x_0)|^2\,dx-\int_{B_{\frac{1}{2}}}\sigma_{\tilde{E}_h}|\D \tilde{u}_hA^{-\frac{1}{2}}(x_0)|^2\,dx\\
& +c\big(n,\alpha,\beta,\lambda,\Lambda,\kappa,\mu,\norm{\D u}_{L^2(\Omega)}\big)\bigg(\tilde \kappa\tilde{r}_h^\mu+[A]_{C^\mu}\tilde{r}_h^\mu\bigg)\frac{1}{2^{n-1+\mu}}.
\end{align}
Using the first variation formula for the perimeter and Proposition \ref{VarB}, for $h$ sufficiently large, we get:
\begin{align}\label{FV1}
P(\tilde{E}_h;B_{\frac{1}{2}})-P(\Phi_{h}(\tilde{E}_h);B_{\frac{1}{2}})=\big(\tilde{r}_h^{\omega\mu}+O\big(\tilde{r}_h^{2\omega\mu}\big)\big)\int_{\dd\tilde{E}_h\cap B_{\frac{1}{2}}}\prodscal{\nu_{\tilde{E}_h}}{e_n}\prodscal{\D'\phi}{\nu'_{\tilde{E}_h}}\,d\mathcal{H}^{n-1},
\end{align}
and
\begin{align}\label{FV2}
& \int_{B_{\frac{1}{2}}}\sigma_{\Phi_{h}(\tilde{E}_h)}|\D (\tilde{u}_h\circ\Phi^{-1}_{h})A^{-\frac{1}{2}}(x_0)|^2\,dx-\int_{B_{\frac{1}{2}}}\sigma_{\tilde{E}_h}|\D \tilde{u}_hA^{-\frac{1}{2}}(x_0)|^2\,dx\notag\\
& \leq c\big(\tilde{r}_h^{\omega\mu}+o\big(\tilde{r}_h^{\omega\mu}\big)\big)\int_{B_{\frac{1}{2}}}|\D \tilde{u}_h|^2\,dx,
\end{align}
for some $c=c(\beta,\lambda,\D\phi,\D\psi)>0$. Inserting \eqref{FV2} and \eqref{FV1} in \eqref{FV0}, dividing by $\sqrt{\varepsilon_h}\big(\tilde{r}_h^{\omega\mu}+O\big(\tilde{r}_h^{2\omega\mu}\big)\big)$ and taking \eqref{EccDirRag} into account, we get
\begin{align*}
&\frac{1}{\sqrt{\varepsilon_h}}\int_{\dd\tilde{E}_h\cap B_{\frac{1}{2}}}\prodscal{\nu_{\tilde{E}_h}}{e_n}\prodscal{\D'\phi}{\nu'_{\tilde{E}_h}}\,d\mathcal{H}^{n-1}\\
& \leq \frac{c}{\sqrt{\varepsilon_h}\big(\tilde{r}_h^{\omega\mu}+O\big(\tilde{r}_h^{\omega\mu}\big)\big)}\bigg(\big(\tilde{r}_h^{\omega\mu}+o\big(\tilde{r}_h^{\omega\mu}\big)\big)\int_{B_{\frac{1}{2}}}|\D \tilde{u}_h|^2\,dx+\tilde{r}_h^{\mu}\bigg)\leq \frac{c}{\sqrt{\varepsilon_h}}\big(\mathcal{D}_{\tilde u_h}(0,1)+\tilde{r}_h^{(1-\omega)\mu}\big)\\
& \leq \frac{c}{\sqrt{\varepsilon_h}}\mathbf{e}(\tilde{E}_h,0,\tilde{\sigma})\leq c\sqrt{\varepsilon_h},
\end{align*}
for some $c=c\big(n,\alpha,\beta,\lambda,\Lambda,\tilde{\kappa},\mu,[A]_{C^\mu},\tilde{\sigma},\tilde{M},\D\phi,\D\psi\big)>0$. Replacing $\phi$ with $-\phi$, we infer that \begin{equation}
\label{5}
\lim_{h\rightarrow +\infty}\frac{1}{\sqrt{\varepsilon_h}} \bigg|\int_{\dd \tilde{E}_h\cap B_{\frac{1}{2}}}\langle \nu_{\tilde{E}_h}, e_n\rangle \langle \D'\phi,\nu_{\tilde{E}_h}' \rangle\,d\mathcal{H}^{n-1}\bigg|=0.
\end{equation}
Decomposing $\dd \tilde{E}_h\cap B_{\frac{1}{2}}=\big([\Gamma_{f_h}\cup(\dd \tilde{E}_h\setminus\Gamma_{f_h})]\setminus(\Gamma_{f_h}\setminus\dd \tilde{E}_h)\big)\cap B_{\frac{1}{2}}$, we deduce
\begin{equation}
\label{4}
\begin{split}
& -\frac{1}{\sqrt{\varepsilon_h}}\int_{\dd \tilde{E}_h\cap B_{\frac{1}{2}}}\langle \nu_{\tilde{E}_h}, e_n\rangle \langle \D'\phi,\nu_{\tilde{E}_h}' \rangle\,d\mathcal{H}^{n-1}
 =\frac{1}{\sqrt{\varepsilon_h}}\bigg[-\int_{ \Gamma_{f_h}\cap B_{\frac{1}{2}}}\langle \nu_{\tilde{E}_h}, e_n\rangle \langle \D'\phi,\nu_{\tilde{E}_h}' \rangle\,d\mathcal{H}^{n-1}\\
& - \int_{(\dd \tilde{E}_h\setminus \Gamma_{f_h})\cap B_{\frac{1}{2}}}\langle \nu_{\tilde{E}_h}, e_n\rangle \langle \D'\phi,\nu_{\tilde{E}_h}' \rangle\,d\mathcal{H}^{n-1}+\int_{(\Gamma_{f_h}\setminus \dd \tilde{E}_h)\cap B_{\frac{1}{2}}}\langle \nu_{\tilde{E}_h}, e_n\rangle \langle \D'\phi,\nu_{\tilde{E}_h}' \rangle\,d\mathcal{H}^{n-1}\bigg].
\end{split}
\end{equation}
Since by the second inequality in \eqref{1} we have
\begin{equation*}
\bigg |\frac{1}{\sqrt{\varepsilon_h}} \int_{(\dd \tilde{E}_h\setminus \Gamma_{f_h})\cap B_{\frac{1}{2}}}\langle \nu_{\tilde{E}_h}, e_n\rangle \langle \D'\phi,\nu_{\tilde{E}_h}' \rangle\,d\mathcal{H}^{n-1} \bigg | \leq c_8\sqrt{\varepsilon_h}\sup_{\R^{n-1}}\abs{\D'\phi},
\end{equation*}
\begin{equation*}
\bigg |\frac{1}{\sqrt{\varepsilon_h}} \int_{(\Gamma_{f_h}\setminus \dd \tilde{E}_h)\cap B_{\frac{1}{2}}}\langle \nu_{\tilde{E}_h}, e_n\rangle \langle \D'\phi,\nu_{\tilde{E}_h}' \rangle\,d\mathcal{H}^{n-1} \bigg | \leq c_8\sqrt{\varepsilon_h}\sup_{\R^{n-1}}\abs{\D'\phi},
\end{equation*}
then, by \eqref{5} and the area formula, we infer
\begin{equation*}
0=\lim_{h\rightarrow +\infty}\frac{-1}{\sqrt{\varepsilon_h}}\int_{ \Gamma_{f_h}\cap B_{\frac{1}{2}}}\langle \nu_{\tilde{E}_h}, e_n\rangle \langle \D'\phi,\nu_{\tilde{E}_h}' \rangle\,d\mathcal{H}^{n-1}=\lim_{h\rightarrow +\infty}\frac{1}{\sqrt{\varepsilon_h}}\int_{ {\mathbf{D}_{\frac{1}{2}}}}\frac{\langle \D' f_h,\D'\phi\rangle}{\sqrt{1+\abs{\D' f_h}^2}}\,dx'.
\end{equation*}
This proves that $g$ is harmonic in $ {\mathbf{D}_{\frac{1}{2}}}$.\\
\indent \textbf{Step 2.} The proof of this step now follows exactly as in \cite{FJ} using the height bound lemma and the reverse Poincaré inequality. We give here the proof for the sake of completeness. Setting
\begin{equation*}
b_h:=\frac{(f_h)_{4\tilde{\tau}}}{\sqrt{1+\abs{(\D' f_h)_{4\tilde{\tau}}}^2}}, \quad \nu_h:=\frac{(-(\D' f_h)_{4\tilde{\tau}},1)}{\sqrt{1+\abs{(\D' f_h)_{4\tilde{\tau}}}^2}}.
\end{equation*}
We want to estimate from above the flatness of $\dd\tilde{E}_h$ towards the hyperplane $\{y\in\R^n\,:\,\prodscal{y}{\nu_h}=b_h\}$ in $B_{4\tilde{\tau}}$ with the excess. More precisely, we show that
\begin{equation}\label{E1}
\limsup_{h\rightarrow +\infty}\frac{1}{\varepsilon_h\tilde{\tau}^{n+1}}\int_{\dd \tilde{E}_h\cap B_{4\tilde{\tau}}}\abs{\langle\nu_h,x\rangle-b_h}^2\,d\mathcal{H}^{n-1}\leq c(n,c_8)\tilde{\tau}^2.
\end{equation}
On one hand, by the mean value property of harmonic functions (see \cite[Lemma 25.1]{Ma}), Jensen's inequality, semicontinuity and the third inequality in \eqref{1} we deduce that
\begin{align}\label{E2}
&\limsup_{h\rightarrow +\infty}\frac{1}{\varepsilon_h\tilde{\tau}^{n+1}}\int_{\dd \tilde{E}_h\cap\Gamma_{f_h}\cap B_{4\tilde{\tau}}}\abs{\langle\nu_h,x\rangle-b_h}^2\,d\mathcal{H}^{n-1}\notag\\
& = \limsup_{h\rightarrow +\infty}\frac{1}{\varepsilon_h\tilde{\tau}^{n+1}}\int_{\dd \tilde{E}_h\cap\Gamma_{f_h}\cap B_{4\tilde{\tau}}}\frac{\abs{\langle -(\D' f_h)_{4\tilde{\tau}},x'\rangle+f_h(x')-(f_h)_{4\tilde{\tau}}}}{1+\abs{(\D' f_h)_{4\tilde{\tau}}}^2}^2\sqrt{1+\abs{\D 'f_h(x')}^2}\,dx'\notag\\
& \leq \limsup_{h\rightarrow +\infty}\frac{1}{\varepsilon_h\tilde{\tau}^{n+1}}\int_{ {\mathbf{D}_{4\tilde{\tau}}}}\abs{f_h(x')-(f_h)_{4\tilde{\tau}}-\langle(\D' f_h)_{4\tilde{\tau}},x'\rangle}^2\,dx'\notag\\
& =\frac{1}{\tilde{\tau}^{n+1}}\int_{ {\mathbf{D}_{4\tilde{\tau}}}}\abs{g(x')-(g)_{4\tilde{\tau}}-\langle(\D' g)_{4\tilde{\tau}},x'\rangle}^2\,dx'\notag\\
& =\frac{1}{\tilde{\tau}^{n+1}}\int_{ {\mathbf{D}_{4\tilde{\tau}}}}\abs{g(x')-g(0)-\langle\D' g(0),x'\rangle}^2\,dx'\notag\\
&\leq c(n)\tilde{\tau}^2\sup_{x'\in  {\mathbf{D}_{4\tilde{\tau}}}}\abs{g(x')-g(0)-\langle\D' g(0),x'\rangle}^2\notag\\
& \leq c(n)\tilde{\tau}^2\int_{ {\mathbf{D}_{\frac{1}{2}}}}\abs{\D' g}^2\,dx'\leq c(n)\tilde{\tau}^2\liminf_{h\rightarrow +\infty}\int_{ {\mathbf{D}_{\frac{1}{2}}}}\abs{\D' g_h}^2\,dx'\leq c(n,c_8)\tilde{\tau}^2,
\end{align}
where we used that $\mathbf{D}_{4\tilde{\tau}}\subset \mathbf{D}_{\frac{1}{4}}$, since $\tilde{\tau}<\frac{1}{16}$. 
On the other hand, from the height bound lemma (see Lemma \ref{height bound}) and \eqref{1}, we immediately get that
\begin{equation}\label{E3}
\lim_{h\rightarrow +\infty}\frac{1}{\varepsilon_h}\int_{(\dd \tilde{E}_h\setminus\Gamma_{f_h})\cap B_{2\tilde{\tau}}}\abs{\langle\nu_h,x\rangle-b_h}^2\,d\mathcal{H}^{n-1}=0.
\end{equation}
Hence, combining \eqref{E2} and \eqref{E3}, we conclude that \eqref{E1} is satisfied.
In order to apply the reverse Poincaré inequality, we show that the sequence $\lbrace {\mathbf e}_h(\tilde{E}_h,4\tilde{\tau},\nu_h) \rbrace_{h\in\N}$ is infinitesimal; indeed, by the definition of excess, Jensen's inequality and the third inequality in \eqref{1} we have
\begin{align*}
& 2(4\tilde{\tau})^{n-1}\limsup_{h\rightarrow+\infty}{\mathbf e}_h(\tilde{E}_h,0,4\tilde{\tau},\nu_h)=\limsup_{h\rightarrow +\infty}\int_{\dd \tilde{E}_h\cap B_{4\tilde{\tau}}}\abs{\nu_{\tilde{E}_h}-\nu_h}^2\,d\mathcal{H}^{n-1}\\
& \leq\limsup_{h\rightarrow +\infty}\bigg[ 2\int_{\dd \tilde{E}_h\cap B_{4\tilde{\tau}}}\abs{\nu_{\tilde{E}_h}-e_n}^2\,d\mathcal{H}^{n-1}+2\abs{e_n-\nu_h}^2\mathcal{H}^{n-1}(\dd \tilde{E}_h\cap B_{4\tilde{\tau}})\bigg]\\
& \leq \limsup_{h\rightarrow +\infty}\bigg[4\varepsilon_h+2\mathcal{H}^{n-1}(\dd \tilde{E}_h\cap B_{4\tilde{\tau}})\frac{\abs{((\D' f_h)_{4\tilde{\tau}},\sqrt{1+\abs{(\D' f_h)_{4\tilde{\tau}}}^2}-1)}^2}{1+\abs{(\D' f_h)_{4\tilde{\tau}}}^2}\bigg]\\
& \leq \limsup_{h\rightarrow +\infty}\big[4\varepsilon_h+4\mathcal{H}^{n-1}(\dd E_h\cap B_{4\tilde{\tau}})\abs{(\D' f_h)_{4\tilde{\tau}}}^2\big]\\
& \leq \limsup_{h\rightarrow +\infty}\bigg[ 4\varepsilon_h+4\int_{ {\mathbf{D}_{\frac{1}{2}}}}\abs{\D' f_h}^2\,dx' \bigg]\leq \lim_{h\rightarrow +\infty}[4\varepsilon_h+4c_8\varepsilon_h]=0.
\end{align*}
Therefore, applying the reverse Poincaré inequality, \eqref{E1} and observing that $\mathbf{C}_{2\tilde{\tau}}\subset B_{4\tilde{\tau}}$, we have for $h$ large that
\begin{align*}
{\mathbf e}(\tilde{E}_h,0,\tau)
& \leq {\mathbf e}(\tilde{E}_h,0,\tau,\nu_h)\\
& \leq c_{11}\bigg(\frac{1}{(2\tilde{\tau})^{n+1}}\int_{\dd\tilde{E}_h\cap \mathbf{C}_{2\tilde{\tau}}(0,\tilde{\tau})}\abs{\langle\nu_h,x\rangle-b_h}^2\,d\mathcal{H}^{n-1}+\mathcal{D}_{\tilde{u}_h}(0,4\tilde{\tau})+\bigl(\tilde{\kappa}+[A]_{C^\mu}\bigr)(2\tilde{\tau} \tilde{r}_h)^\mu\bigg)\\
& \leq \tilde{C}\big(n,\alpha,\beta,\lambda,\Lambda,\kappa,\mu,\norm{\D u}_{L^2(\Omega)}\big)\big(\tilde{\tau}^2{\mathbf e}(\tilde{E}_h,0,1)+\mathcal{D}_{\tilde{u}_h}(0,4\tilde{\tau})+(\tilde{\tau} \tilde{r}_h)^{\mu}\big),
\end{align*}
which is a contradiction if we choose $\overline{C}>\tilde{C}$.

\end{proof}

We use the previous theorem in the proof of the next result.
\begin{Thm}[Excess improvement]\label{MiglEccesso}
For any $\omega\in(0,1)$, $\sigma\in\Big(0,\frac{\lambda^{\frac{1}{2}}}{\Lambda^{\frac{1}{2}}}\Big)$, $M>0$, $\tau\in\Big(0,\frac{\lambda^{\frac{1}{2}}}{16\Lambda^{\frac{1}{2}}}\Big)$ there exists a constant $\varepsilon_6=\varepsilon_6(\sigma,M,\tau)>0$ such that if $(E,u)$ is a $(\kappa,\mu)$-minimizer of $\mathcal{F}_A$ in $B_r(x_0)$, with $x_0\in\dd E$, such that
\begin{equation*}
\mathbf{e}(E,x_0,r)\leq\varepsilon_6,\quad\mathcal{D}_u(x_0,r)+r^{(1-\omega)\mu}\leq M\mathbf{e}(E,x_0,\sigma r),
\end{equation*}
then there exists a constant $c_{14}=c_{14}\big(n,\alpha,\beta,\lambda,\Lambda,\kappa,\mu,\norm{\D u}_{L^2(\Omega)}\big)>0$ such that
\begin{equation*}
\mathbf{e}(E,x_0,\tau r)\leq c_{14}\tau^\mu \big(\mathbf{e}(E,x_0,r)+\mathcal{D}_u(x_0,r)+r^{\mu}\big).
\end{equation*}
\end{Thm}
\begin{proof} Let $\sigma\in\Big(0,\frac{\lambda^{\frac{1}{2}}}{\Lambda^{\frac{1}{2}}}\Big)$, $M>0$, $\tau\in\Big(0,\frac{\lambda^{\frac{1}{2}}}{16\Lambda^{\frac{1}{2}}}\Big)$, and let $(E,u)$ be a $(\kappa,\mu)$-minimizer of $\mathcal{F}_A$ in $B_r(x_0)$, with $x_0\in\dd E$, such that
\begin{equation*}
\mathbf{e}(E,x_0,r)\leq\varepsilon_6,\quad\mathcal{D}_u(x_0,r)+r^{(1-\omega)\mu}\leq M\mathbf{e}(E,x_0,\sigma r).
\end{equation*}
Setting 
\begin{equation}
(\tilde{E},\tilde{u}):=\big(T_{x_0}(E),u\circ T_{x_0}^{-1}\big),\quad \tilde{\kappa}:=\lambda^{-\frac{n}{2}}\kappa,\quad \tilde{r}:=\theta r,
\end{equation}
where $\theta\in\big(0,\min\big\{\Lambda^{-\frac{1}{2}},1\big\}\big)$, by Proposition \ref{Invariance} we have that $(\tilde{E},\tilde{u})$ is a $(\tilde{\kappa},\mu)$-minimizer of $\mathcal{F}_{x_0,A_{x_0}}$ in $B_{\tilde{r}}(x_0)$. By Proposition \ref{CompExcess}, it holds that
\begin{equation*}
\mathbf{e}(\tilde{E},x_0,\tilde{r})\leq \overline{C}_1 \mathbf{e}(E,x_0,r),
\end{equation*}
for some positive constant $\overline{C}_1=\overline{C}_1(n,\lambda,\Lambda)$. Furthermore, estimating
\begin{align*}
\mathcal{D}_{\tilde{u}}(x_0,\tilde{r})
& =\frac{1}{(\theta r)^{n-1}}\int_{B_{\tilde{r}}(x_0)}|\D\tilde{u}|^2\,dy=\frac{\mathrm{det}(A^{-\frac{1}{2}}(x_0))}{(\theta r)^{n-1}}\int_{T^{-1}_{x_0}(B_r(x_0))}|\D u A^{\frac{1}{2}}(x_0)|^2\,dy\\
& \leq\frac{\lambda^{-\frac{n}{2}}\Lambda}{(\theta r)^{n-1}}\int_{B_r(x_0)}|\D u|^2\,dy=\frac{\lambda^{-\frac{n}{2}}\Lambda}{\theta^{n-1}}\mathcal{D}_u(x_0,r)
\end{align*}
and applying again Proposition \ref{CompExcess}, we get
\begin{align*}
& \mathcal{D}_{\tilde{u}}(x_0,\tilde{r})+\tilde{r}^{(1-\omega)\mu}\leq c(n,\lambda,\Lambda)\big(\mathcal{D}_u(x_0,r)+r^{(1-\omega)\mu}\big)\leq c(n,\lambda,\Lambda)M\mathbf{e}(E,x_0,\sigma r)\\
& \leq  c(n,\lambda,\Lambda)M\mathbf{e}\big(\tilde{E},x_0,\sigma\lambda^{-\frac{1}{2}} r\big)\leq \overline{C}_2M\mathbf{e}(\tilde{E},x_0,\tilde{\sigma}\tilde{r}),
\end{align*}
for some positive constant $\overline{C}_2=\overline{C}_2(n,\lambda,\Lambda)$, where $\tilde{\sigma}:=\frac{\lambda^{-\frac{1}{2}}}{\theta}\sigma<1$, since $\sigma<\theta\lambda^{\frac{1}{2}}$. Choosing $\varepsilon_6>0$ such that $\overline{C}_1\varepsilon_6<\tilde{\varepsilon}$ and setting $\tilde{M}:=\overline{C}_2 M$, we apply Theorem \ref{Teorema eccesso raddrizzato} to obtain
\begin{equation*}
\mathbf{e}(\tilde{E},x_0,\tilde{\tau}\tilde{r})\leq C\big(\tilde{\tau}^2\mathbf{e}(\tilde{E},x_0,\tilde{r})+\mathcal{D}_{\tilde{u}}(x_0,4\tilde{\tau}\tilde{r})+(\tilde{\tau}\tilde{r})^{\mu}\big),
\end{equation*}
for some positive constant $C=C\big(n,\alpha,\beta,\lambda,\Lambda,\kappa,\mu,\norm{\D u}_{L^2(\Omega)}\big)$, where $\tilde{\tau}:=\frac{\tau}{\lambda^{\frac{1}{2}}\theta}<\frac{1}{16}$, since $\tau<\frac{\lambda^{\frac{1}{2}}\theta}{16}$. Leveraging Proposition \ref{CompExcess}, we get
\begin{align}
\label{aa1}
\mathbf{e}(E,x_0,\tau r)
& =\mathbf{e}\big(E,x_0,\tilde{\tau}\lambda^{\frac{1}{2}}\tilde{r}\big)\leq c(n,\lambda,\Lambda)\mathbf{e}(\tilde{E},x_0,\tilde{\tau}\tilde{r})\notag\\
& \leq C\big(\tilde{\tau}^2\mathbf{e}(\tilde{E},x_0,\tilde{r})+\mathcal{D}_{\tilde{u}}(x_0,4\tilde{\tau}\tilde{r})+(\tilde{\tau}\tilde{r})^\mu\big).
\end{align}
On one hand, by Proposition \ref{CompExcess} we observe that
\begin{equation}\label{EF1}
\mathbf{e}(\tilde{E},x_0,\tilde{r})\leq c(n,\lambda,\Lambda)\mathbf{e}\big(E,x_0,\Lambda^{\frac{1}{2}}\tilde{r}\big)\leq c(n,\lambda,\Lambda)\mathbf{e}(E,x_0,r),
\end{equation}
being $\Lambda^{\frac{1}{2}}\tilde{r}\leq r$. One the other hand, choosing $\varepsilon_6<\varepsilon_2$, by Proposition \ref{DecayDirichlet} it follows
\begin{equation}\label{EF2}
\mathcal{D}_{\tilde{u}}(x_0,4\tilde{\tau}\tilde{r})\leq c(n,\lambda,\Lambda) \mathcal{D}_u\big(x_0,4\tilde{\tau}\Lambda^{\frac{1}{2}}\tilde{r}\big)\leq c(n,\lambda,\Lambda)\tilde{\tau} \mathcal{D}_u(x_0,r),
\end{equation}
since $4\tilde{\tau}\Lambda^{\frac{1}{2}}\tilde{r}\leq r$, being $\tau<\frac{\lambda^{\frac{1}{2}}}{4\Lambda^{\frac{1}{2}}}$.
Inserting \eqref{EF1} and \eqref{EF2} in \eqref{aa1}, we obtain
\begin{align*}
\mathbf{e}(E,x_0,\tau r)
& \leq C\big(\tilde{\tau}^2\mathbf{e}(E,x_0,r)+\tilde{\tau}\mathcal{D}_u(x_0,r)+\tilde{\tau}^{\mu} r^{\mu}\big)\\
& \leq C\tilde{\tau}^{\mu}\big(\mathbf{e}(E,x_0,r)+\mathcal{D}_u(x_0,r)+ r^{\mu}\big)\\
& \leq C\tau^{\mu}\big(\mathbf{e}(E,x_0,r)+\mathcal{D}_u(x_0,r)+ r^{\mu}\big),
\end{align*}
which is the thesis.
\end{proof}

Leveraging the results proved in the previous sections, we are able to prove Theorem \ref{Teorema principale}.
\begin{proof}[Proof of Theorem \ref{Teorema principale}]
Let $U\Subset\Omega$ be an open set. We prove that for every $\omega\in(0,1)$ and
$\tau\in(0,1)$ there exist two positive constants $\overline{\varepsilon}=\overline{\varepsilon}(\tau,U)$ and $\overline{C}$ such that if $x_0\in\dd E$, $B_r(x_0)\subset U$ and ${\mathbf e}(x_0,r)+\mathcal{D}(x_0,r)+r^{(1-\omega)\mu}<\overline{\varepsilon}$, then
\begin{equation}
\label{Eqn 8}
{\mathbf e}(x_0,\tau r)+\mathcal{D}(x_0,\tau r)+(\tau r)^{(1-\omega)\mu}\leq \overline{C}\tau^{(1-\omega)\mu}\big({\mathbf e}(x_0,r)+\mathcal{D}(x_0,r)+ r^{(1-\omega)\mu}\big).
\end{equation}
We fix $\tau\in(0,1)$. Setting
\begin{equation*}
\overline{\tau}:=\frac{\lambda^{\frac{1}{2}}}{16\Lambda^{\frac{1}{2}}},\quad\overline{\sigma}:=\frac{\lambda^{\frac{1}{2}}}{\Lambda^{\frac{1}{2}}},
\end{equation*}
we may assume without loss of generality that 
\begin{equation*}
\tau<\min\bigg\{\overline{\tau},\frac{\overline{\sigma}}{2}\bigg\}=\overline{\tau}.
\end{equation*}
Furthermore we fix $\sigma:=2\tau<\overline{\sigma}$. We distinguish two cases. \\
\indent\textit{Case 1:} $\mathcal{D}_u(x_0,r)+r^{(1-\omega)\mu}\leq \tau^{-1}{\mathbf e}(x_0,\sigma r)$. Choosing $\overline{\varepsilon}<\varepsilon_6(\sigma,\tau,\tau)$ it follows from Theorem \ref{MiglEccesso} that
\begin{equation*}
{\mathbf e}(x_0,\tau r)\leq c_{14}\tau^{\mu} \big(\mathbf{e}(x_0,r)+\mathcal{D}_u(x_0,r)+r^{\mu}\big).
\end{equation*}
Furthermore, choosing $\overline{\varepsilon}<\varepsilon_2(\tau)$, applying Proposition \ref{DecayDirichlet} we get \eqref{Eqn 8}.

\textit{Case 2:} ${\mathbf e}(x_0,\sigma r)\leq\tau\big(\mathcal{D}_u(x_0, r)+r^{(1-\omega)\mu}\big)$. By the property of the excess at different scales, we infer
\begin{equation*}
{\mathbf e}(x_0,\tau r)\leq 2^{n-1}{\mathbf e}(x_0,\sigma r)\leq 2^{n-1}\tau\big(\mathcal{D}_u(x_0,r)+r^{(1-\omega)\mu}\big),
\end{equation*}
obtaining \eqref{Eqn 8}.\\
\indent Thus, choosing $\overline{\varepsilon}=\min\{ \varepsilon_2(\tau),\varepsilon_6(2\tau,\tau,\tau)\}$, we conclude that the inequality \eqref{Eqn 8} is verified.\\
 We fix $\sigma\in\big(0,\frac{(1-\omega)\mu}{2}\big)$ and choose $\tau_0\in(0,1)$ such that $\overline{C}\tau_0^{(1-\omega)\mu}\leq\tau_0^{2\sigma}$ and we define
\begin{equation*}
\Gamma\cap U:=\big\{ x\in\dd E\cap U\,:\, {\mathbf e}(x,r)+\mathcal{D}(x,r)+r^{(1-\omega)\mu}<\overline{\varepsilon}(\tau_0,U),
\text{ for some }r>0\text{ such that }B_r(x_0)\subset U \big\}.
\end{equation*}
We note that $\Gamma\cap U$ is relatively open in $\dd E$. We show that $\Gamma\cap U$ is a $C^{1,\sigma}$-hypersurface. Indeed, inequality \eqref{Eqn 8} implies via standard iteration argument that if $x_0\in\Gamma\cap U$ there exist $r_0>0$ and a neighborhood $V$ of $x_0$ such that for every $x\in\dd E\cap V$ it holds:
\begin{equation*}
{\mathbf e}\big(x,\tau_0^k r_0\big)+\mathcal{D}\big(x,\tau_0^k r_0\big)+\big(\tau_0^k r_0\big)^{(1-\omega)\mu}\leq \tau_0^{2\sigma k}, \quad\text{for }k\in\N_0.
\end{equation*}
In particular, ${\mathbf e}(x,\tau_0^k r_0)\leq \tau_0^{2\sigma k}$ and, arguing as in \cite{FJ}, we obtain that for every $x\in\dd E\cap V$ and $0<s<t<r_0$ it holds
\begin{equation*}
|(\nu_E)_s(x)-(\nu_E)_t(x)|\leq ct^\sigma,
\end{equation*}
for some constant $c=c(n,\tau_0,r_0)$, where
\begin{equation*}
(\nu_E)_t(x)=\dashint_{\dd E\cap B_t(x)}\nu_E\,d\mathcal{H}^{n-1}.
\end{equation*}
The previous estimate first implies that $\Gamma\cap U$ is $C^1$. By a standard argument we then deduce again from the same estimate that $\Gamma\cap U$ is a $C^{1,\sigma}$-hypersurface. Since $\omega$ is arbitrary, we gain that $\Gamma$ is a $C^{1,\sigma}$-hypersurface, for any $\sigma\in (0,\frac{\mu}{2})$. We define $\Gamma:=\cup_i(\Gamma\cap U_i)$, where $(U_i)_i$ is an increasing sequence of open sets such that $U_i\Subset\Omega$ and $\Omega=\cup_i U_i$. We are left to prove that there exists $\varepsilon>0$ such that
\begin{equation*}
\mathcal{H}^{n-1-\varepsilon}(\dd E\setminus \Gamma)=0.
\end{equation*}
Setting $\Sigma=\Big\{ x\in\dd E\setminus\Gamma\,:\, \lim\limits_{r\rightarrow 0}\mathcal{D}(x,r)=0 \Big\}$, by \cite[Lemma 2.1]{FJ} we have that $\D u\in L^{2s}_{loc}(\Omega)$ for some $s=s(n,\alpha,\beta)>1$ and we have that
\begin{equation*}
\text{dim}_{\mathcal{H}}\Big(\Big\{x\in\Omega\,:\, \limsup_{r\rightarrow 0}\mathcal{D}(x,r)>0 \Big\}\Big)\leq n-s.
\end{equation*}
The conclusion follows in a standard way as in \cite{FJ} (see also \cite{DFR} and \cite{DF0}) showing that $\Sigma=\emptyset$ if $n\leq 7$ and $\text{dim}_\mathcal{H}(\Sigma)\leq n-8$ if $n\geq 8$. In both cases, Lemma \ref{Lemma compattezza} will be employed.
\end{proof}

\section{An application to a costrained problem}
\label{From costrained to penalized problem}
In this section we show an application of Theorem \ref{Teorema principale} to the following costrained problem:
\begin{equation}
\label{P_c}
\min_{\substack{E\in \mathcal{A}(\Omega)\\v\in u_0+H_0^1(\Omega)}}
\left\{
\mathcal{F}_A(E,v;\Omega)\,:\,  |E|=d
\right\},
\tag{$P_c$}
\end{equation}
where $u_0\in H^1(\Omega)$, $0<d<|\Omega|$ are given and $\mathcal{A}(\Omega)$ is the class of all subsets of $\Omega$ with finite perimeter. We assume that $\Omega$ is connected.\\
\indent In this perspective we need to distinguish the H\"older exponent of the matrix $A$, which we denote here with $\gamma$, from the exponent $\mu$ appearing in the Definition \ref{k-mu}. In Theorem \ref{Teorema Penalizzazione}, we show that, for sufficiently large values of $\kappa>0$, minimizing couples of \eqref{P_c} are solutions of the following penalized problem
\begin{equation}
\label{P}
\min_{\substack{E\in \mathcal{A}(\Omega)\\v\in u_0+H_0^1(\Omega)}}
\mathcal{F}_{\kappa}(E,v;\Omega),
\tag{$P$}
\end{equation}
where the functional $\mathcal{F}_\kappa$ is defined by
\begin{equation}\label{intro1}
{\mathcal F}_{\kappa}(E,v;\Omega):=\mathcal{F}_A(E,v;\Omega)+\kappa\big||E| -d\big|^{\gamma}.
\end{equation}
Now we prove the penalization theorem. For simplicity of notation, we denote
\begin{equation*}
a(x,\nu)=\prodscal{A(x)\nu}{\nu}^{\frac{1}{2}},\quad \forall x,\nu\in \R^n.
\end{equation*}
It will be advantageous to have some estimates about the dependence of the integrand $a$ on $x$ and $\nu$.
\begin{Rem}[Continuity of $a$ with respect to $x$ and $\nu$] It is straightforward to check that the following inequalities hold:
\begin{align}
\label{Depx}& |a(x,\nu)-a(y,\nu)|\leq \frac{1}{2\sqrt{\lambda}}[A]_{C^{\mu}}|x-y|^{\mu}, \quad \forall x,y\in \R^n, \quad |\nu|=1.\\
\label{Depv}& |a(x,\xi)-a(x,\eta)|\leq \frac{\Lambda}{\sqrt{\lambda}}|\xi-\eta|,\quad \forall x\in \R^n, \quad\forall\xi,\eta\in\R^n.
\end{align}
\end{Rem}

 As explained before, the proof of the equivalence between the solution of the constrained problem $\eqref{P_c}$  and the penalized problem $\eqref{P}$ follows. We adapt a result proved in \cite{EF} to our setting.
\begin{Thm}
\label{Teorema Penalizzazione} There exists $\kappa_{0}>0$ such that if $(E,u)$ is a minimizer of the functional
\begin{equation}
\label{Penalized}
{\mathcal F}_{\kappa}(A,w)=\int_\Om\sigma_A|\D u|^2\,dx +\mathbf{\Phi}_A(E;\Omega)+\kappa\big||A| - d\big|^\gamma,
\end{equation}
for some $\kappa \geq \kappa_0>0$, among all configurations $(A,w)$ such that $w=u_0$ on $\partial \Omega$,
then $|E|=d$ and $(E,u)$ is a minimizer of problem $\eqref{P_c}$.
Conversely, if $(E,u)$ is a minimizer of problem $\eqref{P_c}$, then it is a minimizer of \eqref{Penalized}, for all $\kappa \geq \kappa_0$.
\end{Thm}
\begin{proof}
The argument is very similar to the one in \cite[Theorem 1]{EF} (see also \cite{EL1} and \cite{EL2}). For reader's convenience, we give here the sketch of the proof, emphasizing main ideas and some differences with respect to the case treated in \cite{EF}.\\
\indent The first part of the theorem can be proved by contradiction. We assume that there exist a positive sequence $(\kappa_h)_{h\in \mathbb N}$ such that $\kappa_h\rightarrow +\infty$, as $h\rightarrow +\infty$, and a sequence of configurations $(E_h,u_h)$ minimizing $\mathcal{F}_{\kappa_h}$ and such that $u_h=u_0$ on $\partial \Omega$ and $|E_h|\neq d$, for all $h\in\N$. We choose an arbitrary fixed $E_0\subset \Omega$ with finite perimeter and such that $|E_0|=d$. We point out that 
\begin{equation}\label{Theta}
\mathcal{F}_{\kappa_h}(E_h,u_h)\leq\mathcal{F}(E_0,u_0):=\Theta.
\end{equation}
Without loss of generality we can assume that $|E_h|<d$, the case $|E_h|>d$ being similar. Our aim is to show that for $h$ sufficiently large, there exists a configuration $(\widetilde{E}_h,\tilde{u}_h)$ such that $\mathcal{F}_{\kappa_h}(\widetilde{E}_h,\tilde{u}_h)< \mathcal{F}_{\kappa_h}({E_h},{u_h})$, thus proving the result by contradiction.\\
\indent By condition $\eqref{Theta}$, it follows that the sequence $(u_h)_{h\in\N}$ is bounded in $H^1(\Omega)$, the perimeters of the sets $E_h$ in $\Omega$ are unifromly bounded and $|E_h|\rightarrow d$. Therefore, possibly extracting a not relabelled subsequence, we may assume that there exists a configuration $(E,u)$ such that $u_h\rightarrow u$ weakly in $H^1(\Omega)$, $\mathbbm{1}_{E_h}\rightarrow \mathbbm{1}_{E}$ a.e. in $\Omega$, where the set $E$ is of finite perimeter in $\Omega$ and $|E|=d$. The couple $(E,u)$ will be used as a reference configuration for the definition of $(\widetilde{E}_h,\tilde{u}_h)$.\\

{\bf Step 1.} {\em Construction of $(\widetilde{E}_h,\tilde{u}_h)$}.
Proceeding exactly as in \cite{EF}, since $\Omega$ is connected, we can take a point $x\in\partial^*E\cap\Om$. We oobserve that,  
given $\varepsilon> 0$ sufficiently small, we can find around $x$ a point
$x'$ and $r > 0$ such that


$$
|E\cap B_{r/2}(x')|<\e r^n,\qquad|E\cap B_r(x')|>\frac{\omega_nr^n}{2^{n+2}}\,.
$$
We assume without loss of generality that $x'=0$, and from now on we denote by $B_r$ the balls centered at the origin. From the convergence of $E_h$ to $E$ we have that, for $h$ sufficiently large,
\begin{equation}\label{unodue}
|E_h\cap B_{r/2}|<\e r^n,\qquad|E_h\cap B_r|>\frac{\omega_nr^n}{2^{n+2}}\,.
\end{equation}
Now we define the following bi-Lipschitz map used in \cite{EF} which maps $B_r$ into itself: 
\begin{equation}\label{unotre}
f(x):=
\begin{cases}
\bigl(1-\sigma_h\big(2^n-1\big)\bigr)x & \text{if}\,\,\,|x|<\displaystyle\frac{r}{2},\\
x+\sigma_h\bigg(\displaystyle 1-\frac{r^n}{|x|^n}\bigg)x & \text{if}\,\,\,\displaystyle\frac{r}{2}\leq|x|<r,\\
x & \text{if}\,\,\,|x|\geq r\,,
\end{cases}
\end{equation}
for some $0<\sigma_h<1/2^n$ such that, setting
$$
{\widetilde E}_h=f(E_h),\qquad{\tilde u}_h=u_h\circ f^{-1}\,,
$$
we have $|{\widetilde E}_h|<d$. By minimality, We obtain
\begin{align}\label{unoquattro}
{\mathcal F}_{\kappa_h}(u_h,E_h)-{\mathcal F}_{\kappa_h}({\tilde u}_h,{\widetilde E}_h) &=\biggl[\int_{B_r}\sigma_{E_h}|\D u_h|^2\,dx-\int_{B_r}\sigma_{\tilde{E}_h}|\D \tilde{u}_h|^2\,dx\biggr] \nonumber
\\ & +\bigl[\mathbf{\Phi}_A(E_h,{\overline B}_r)-\mathbf{\Phi}_A({\widetilde E}_h,{\overline B}_r)\bigr]+\kappa_h\bigl[(d-|{ E}_h|)^\gamma-(d-|\widetilde{E}_h|)^\gamma\bigr]\nonumber \\
&= I_{1,h}+I_{2,h}+I_{3,h}. 
\end{align}
For simplicity of notation we will denote in the following 
\begin{equation*}
    g(y)=f^{-1}(y),\quad\forall y\in \R^n.
\end{equation*} 
We will use in the sequel some estimates 
for the map $f$ that can be easily obtained by direct computation (see \cite{EF} for the explicit calculation). These estimates are trivial for $|x|<r/2$, whereas they
can be deduced by the explicit expression of $\nabla f$ for $r/2<|x|<r$ , that is
\begin{equation*}
\frac{\partial f_i}{\partial x_j}(x)=\delta_{ij}+\sigma_h\Bigl[\bigl(1-\frac{ r^n}{|x|^n}\bigr)\delta_{ij}+nr^n\frac{x_ix_j}{|x|^{n+2}}\Bigr],\quad \forall i,j\in\{1,\dots, n\}.
\end{equation*}
There exists a constant $C=C(n)$ depending only on $n$ such that, 
\begin{equation}\label{phi}
\bigl\|\nabla g-I\bigr\|\leq C(n)\sigma_h, \quad\forall y\in B_r,
\end{equation}
\begin{equation}
\label{Jphi}
1+C(n)\sigma_h\leq Jf(x)\leq1+2^n n\sigma_h,\quad\forall x\in B_r.
\end{equation}

{\bf Step 2.} {\em Estimate of $I_{1,h}$}. 
Performing the change of variables $y=f(x)$, and observing that $\mathbbm{1}_{\widetilde{E}_h}\circ f=\mathbbm{1}_{E_h}$, we get
\begin{equation}
I_{1,h}=\int_{B_r}\sigma_{E_h}(x)\bigl[|\nabla u_h(x)|^2-\bigl|\nabla u_h(x)\circ\nabla f^{-1}(f(x))\bigl|^2Jf(x)\bigr]\,dx.
\end{equation}
By means of the same computation as in \cite{EF}, Using \eqref{phi} and \eqref{Jphi} we deduce that
\begin{equation}
\label{I1h}
I_{1,h}\geq - \overline{C}_1\Theta\sigma_h,
\end{equation}
for some positive constant $\overline{C}_1=\overline{C}_1(n)$.\\
\indent {\bf Step 3.} {\em Estimate of $I_{2,h}$}. In order to estimate $I_{2,h}$ we can use a generalized area formula for maps between rectifiable sets involving anistropies.
We recall that (see \cite[Proposition 17.1]{Ma}), if $E$ is a set of locally finite
perimeter in $\R^n$, then $f(E)$ is a set of locally finite perimeter in $\R^n$ and 
\begin{equation}\label{chnormal}
\partial^{*}f(E)=f(\partial^*E),\quad\quad\nu_{f(E)}(y)=\frac{[\nabla g(y)]^t(\nu_E(y)}{|[\nabla g(y)]^t(\nu_E(y))|},\quad\forall y\in \partial^{*}f(E).
\end{equation}
Using \cite[formula (17.6)]{Ma} we can easily deduce that
\begin{equation}\label{Diffeo1}
\int_{\partial^*f(E)}\phi(f^{-1}(y))\, d\mathcal{H}^{n-1}_y=\int_{\partial^*E}\phi(x)Jf(x)|(\nabla g\circ f)^t\nu_E(x)|\, d\mathcal{H}^{n-1}_x
\end{equation}
for any Borel function $\phi$ defined on $\partial^*E$. If we choose $\phi(x)=a(f(x), \nu_{f(E)}(f(x)))$ in \eqref{Diffeo1} we deduce that
\begin{equation}\label{Diffeo2}
\mathbf{\Phi}_A(f(E))=\int_{\partial^*E} a\bigl(f(x),[\nabla g(f(x))]^t\nu_E(x)\bigr)Jf(x)
\, d\mathcal{H}^{n-1}_x.
\end{equation}
For the proof of the aforementioned formula in a more general framework the reader is addressed to \cite[Proposition A.1]{Simm}.
Now we are ready to estimate the following quantity:
\begin{align}
I_{2,h}
& =\bigl[\mathbf{\Phi}_A(E_h,{\overline B}_r)-\mathbf{\Phi}_A({\widetilde E}_h,{\overline B}_r)\bigr]\\
& =\int_{\partial^*E_h\cap{\overline B}_r}\bigl[ a\bigl(f(x),[\nabla g(f(x))]^t\nu_{E_h}(x)\bigr)-a\bigl(f(x),\nu_{E_h}(x)\bigr)\bigr]Jf(x)
\, d\mathcal{H}^{n-1}_x\\
& +\int_{\partial^*E_h\cap{\overline B}_r}\bigl[ a\bigl(f(x),\nu_{E_h}(x)\bigr)-a\bigl(x,\nu_{E_h}(x)\bigr)\bigr]Jf(x)
\, d\mathcal{H}^{n-1}_x\\
& +\int_{\partial^*E_h\cap{\overline B}_r}\bigl[Jf(x)-1\bigr] a\bigl(x,\nu_{E_h}(x)\bigr)
\, d\mathcal{H}^{n-1}_x=J_{1,h}+J_{2,h}+J_{3,h}.
\end{align}
Using \eqref{phi} and \eqref{Depv} we deduce
\begin{equation}\label{J1}
|J_{1,h}|\leq \frac{\Lambda}{\sqrt{\lambda}}\int_{\partial^*E_h\cap{\overline B}_r}Jf(x)\bigl|\bigl[(\nabla g)^t-I\bigr]\nu_{E_h}(x)\bigr|\, d\mathcal{H}^{n-1}_x\leq C(n)\Theta\frac{\Lambda}{\sqrt{\lambda}}\sigma_.
\end{equation}
Applying \eqref{Depx} we obtain
\begin{equation}\label{J2}
|J_{2,h}|\leq \frac{[A]_{C^{\gamma}}}{2\sqrt{\lambda}}\int_{\partial^*E_h\cap{\overline B}_r}Jf(x)|f(x)-x|^{\gamma}\
d\mathcal{H}^{n-1}_x\leq \frac{[A]_{C^{\gamma}}}{2\sqrt{\lambda}} C(n)\Theta\sigma^{\gamma}.
\end{equation}
Finally from \eqref{Jphi} we have
\begin{equation}
|J_{3,h}|\leq \int_{\partial^*E_h\cap{\overline B}_r}n2^n a(x,\nu_{E_h}(x))\sigma_h\
d\mathcal{H}^{n-1}_x\leq n2^n\Theta\sqrt{\Lambda}\sigma_h.
\end{equation}
Summarizing we conclude that
\begin{equation}\label{I2h}
I_{2,h}\geq-\overline{C}_2\Theta \sigma_h^{\gamma},
\end{equation}
for some positive constant $\overline{C}_3=\overline{C}_3\big(n,\lambda,\Lambda,[A]_{C^{\gamma}}\big)$.\\
\indent {\bf Step 4.} {\em Estimate of $I_{3,h}$}.
The following estimate is contained in \cite[Theorem 2]{EL2} and we detail it for reader's convenience.\\
\indent First we recall \eqref{unodue}, \eqref{unotre}, \eqref{Jphi}, thus getting
\begin{align*}
|\tilde{E}_h|-|E_h|
& = \int_{E_h\cap B_r\setminus B_{r/2}}\!\left(J f(x)-1\right)\,dx+\int_{E_h\cap B_{r/2}}\!\left(J f(x)-1\right)\,dx\\
& \geq\bigg(\frac{\omega_n}{2^{n+2}}-\e\bigg)\sigma_h r^n-\bigl[1-\bigl(1-(2^n-1)\sigma_h\bigr)^n\bigr]\e r^n \\
& \geq \sigma_h r^n\bigg[\frac{\omega_n}{2^{n+2}}-\e-(2^n-1)n\e\bigg].
\end{align*}
Therefore, if we choose $0<\e<\overline{\varepsilon}(n)$, for some $\overline{\varepsilon}$ sufficiently small, we have that
\begin{equation}\label{I3hh}
\kappa_h(|\tilde{E}_h|-|E_h|)\geq\lambda_h C(n)\sigma_h r^n.
\end{equation}
Moreover, if denoting $\delta_h:=d-|E_h|$, we choose $\sigma_h$ in such a way that
$|\tilde{E}_h|-|E_h|\leq \delta_h/2 $ thus respecting the condition $|{\widetilde E}_h|<d$.
Taking this into account, proceding as before and using \eqref{Jphi}, we have
\begin{align*}
|\tilde{E}_h|-|E_h|
& = \int_{E_h\cap B_{r}}\!\left(J f(x)-1\right)\,dx
\leq n2^n \sigma_h r^n.
\end{align*}
Then we choose $\sigma_h$ such that
\begin{equation}
\delta_h\leq \sigma_h\leq \frac{\delta_h}{n2^{n+1}r^n}.
\end{equation}
We remark that in the last condition we imposed also that $\sigma_h$ is comparable with $\delta_h$, which is crucial in the following estimate.
Resuming \eqref{I3hh} we can conclude
\begin{align}\label{I3h}
I_{3,h}
& =\kappa_h\bigl[(d-|{ E}_h|)^\gamma-(d-|\widetilde{E}_h|)^\gamma\bigr]\geq \kappa_h\frac{\gamma}{(d-|E_h|)^{1-\gamma}}(|\tilde{E}_h|-|E_h|)\\
& = \kappa_h\gamma(d-|E_h|)^{\gamma}\frac{|\tilde{E}_h|-|E_h|}{d-|E_h|}\geq \kappa_h\gamma \delta_h^\gamma\frac{c_2(n)\sigma_h r^n}{\delta_h}\\
& \geq \kappa_h \overline{C}_3\sigma_h^\gamma r^n,
\end{align}
for some positive constant $\overline{C}_3=\overline{C}_3(n,\gamma)$.\\

 From the previous inequality inequality, recalling \eqref{unoquattro}, \eqref{I1h} and \eqref{I2h}, we obtain
\begin{equation}
{\mathcal F}_{\kappa_h}(u_h,E_h)-{\mathcal F}_{\kappa_h}({\tilde u}_h,{\widetilde E}_h)\geq\sigma_h^{\gamma}\bigl(\kappa_h\overline{c_3}r^n-\Theta(\overline{C}_1+\overline{C}_2)\bigr)>0,
\end{equation}
if $\kappa_h$ is sufficiently large. This contradicts the minimality of $(E_h,u_h)$, thus concluding the proof.
\end{proof}
\begin{Rem}
Theorem \ref{Teorema Penalizzazione} allows us to prove the regularity of solutions of the free boundary problem under the constraint $|E|=d$. Under the assumption
\begin{equation}
    \gamma\in\bigg(\frac{n-1}{n},1\bigg),
\end{equation}
the parameter $\mu:=\gamma n-n+1$ is positive and, by Theorem \ref{Teorema Penalizzazione}, any minimizing couple $(E,u)$ of \eqref{P_c} is a $(\kappa,\mu)$-minimizer of $\mathcal{F}_\kappa$, for $\kappa\geq\kappa_0$, where $\kappa_0$ is the constant appearing in Theorem \ref{Teorema Penalizzazione}. Thus, we are in position to implement the regularity theory of the previous sections to $(E,u)$ by applying Theorem \ref{Teorema principale}.
\end{Rem}

\emph{Acknowledgements} The authors are members of the Gruppo Nazionale per l’Analisi Matematica, la Probabilità e le loro Applicazioni (GNAMPA) of the Istituto Nazionale di Alta Matematica (INdAM). The authors wish to acknowledge financial support from IN-dAM GNAMPA Project 2024 "Regolarità per problemi a frontiera libera e disuguaglianze funzionali in contesto finsleriano". The work of GP was partially supported by the project "Start" within the program of the University "Luigi Vanvitelli" reserved to young researchers, Piano strategico 2021-2023\\

\textbf{Data Availability} Data sharing not applicable to this article as no datasets were generated or analysed during the current study.\\

\textbf{Declarations}\\

\textbf{Conflicts of interest} The authors declare that they have no conflicts of interest.

\noindent
\author{Luca Esposito},
Dipartimento di Matematica, Università degli Studi di Salerno, Via Giovanni Paolo II 132, Fisciano 84084, Italy\\
luesposi@unisa.it

\medskip
\noindent
\author{Lorenzo Lamberti},
Dipartimento di Matematica, Università degli Studi di Salerno, Via Giovanni Paolo II 132, Fisciano 84084, Italy\\
llamberti@unisa.it

\medskip
\noindent
\author{Giovanni Pisante}, Dipartimento di Matematica e Fisica, Universit\`{a} della Campania {\em Luigi Vanvitelli},
viale Lincoln 5, Caserta 81100, Italy\\
\noindent
giovanni.pisante@unicampania.it

\end{document}